\documentclass[preprint,3p]{elsarticle}
\usepackage{stmaryrd}
\usepackage{amsmath}
\usepackage{amssymb}
\usepackage{tabularx}

\usepackage{graphicx}
\usepackage{epstopdf}
\usepackage{amssymb}
\usepackage{color}
\usepackage{amsmath, amsthm}

\newcommand{\p}{\partial}

\newcommand{\Og}{\Omega}
\newcommand{\fl}[2]{\frac{#1}{#2}}
\newcommand{\dt}{\delta}
\newcommand{\nn}{\nonumber}
\newcommand{\ap}{\alpha}
\newcommand{\bt}{\beta}

\newcommand{\Dt}{\Delta}

\newcommand{\be}{\begin{equation}}
\newcommand{\ee}{\end{equation}}
\newcommand{\ba}{\begin{array}}
\newcommand{\ea}{\end{array}}
\newcommand{\bea}{\begin{eqnarray}}
\newcommand{\eea}{\end{eqnarray}}
\newcommand{\beas}{\begin{eqnarray*}}
\newcommand{\eeas}{\end{eqnarray*}}

\newtheorem{remark}{Remark}[section]
\newcommand{\bx}{{\bf x} }
\newcommand{\tbx}{\widetilde{\bf x}}
\newcommand{\tx}{\widetilde{x}}
\newcommand{\ty}{\widetilde{y}}
\newcommand{\tz}{\widetilde{z}}
\newcommand{\gm}{{\gamma} }

\title{An efficient spectral method for computing dynamics of rotating two-component Bose--Einstein
condensates via coordinate transformation}

\author[1]{Ming Ju}
\ead{jming@csrc.ac.cn}
\address[1]{Beijing Computational Science Research Center, No. 3 He-Qing Road, Hai-Dian
District, Beijing, P.R. China 100084}

\author[2]{Qinglin Tang}
\ead{g0800880@nus.edu.sg}
\address[2]{Department of Mathematics and Center for Computational Science and Engineering,
National University of Singapore, Singapore 119076}

\author[3]{Yanzhi Zhang\corref{cor}}
\ead{zhangyanz@mst.edu}
\cortext[cor]{Corresponding author}
\address[3]{Department of Mathematics and Statistics, Missouri University of Science
and Technology, Rolla, MO 65409-0020, USA}

\begin{document}

\begin{abstract}
In this paper,  we propose an efficient and accurate numerical method for computing the dynamics
of rotating two-component Bose--Einstein condensates (BECs)  which is described by  coupled
Gross--Pitaevskii equations (CGPEs) with an angular momentum rotation term and an external driving
field.  By introducing rotating Lagrangian coordinates, we eliminate the angular momentum 
rotation term  from the CGPEs, which allows us to develop an efficient numerical method. Our 
method has  spectral accuracy in all spatial dimensions and moreover it can be easily implemented in 
practice. 
To examine its performance, we compare our method with those reported in literature.  
Numerical results show that to achieve the same accuracy, our method needs much shorter 
computing time.  We also applied our method to study the dynamic properties of rotating 
two-component BECs.  Furthermore,  we generalize our method to solve the vector Gross--Pitaevskii 
equations (VGPEs)   which is used to study rotating multi-component BECs.
\end{abstract}


\maketitle
\section{Introduction}
\label{section1}


The Bose--Einstein condensation (BEC),  which affords an  astonishing  glimpse into the macroscopic 
quantum world, has been extensively studied 
since its first realization in 1995 \cite{Anderson1995, Bradley, Davis1995}.
Later, with the observation of quantized vortices in BECs \cite{Matthews1999, Madison2000}, 
attention has been broaden to explore vortex states and their dynamics  associated with superfluidity.  
Rotating BECs which are known to exhibit highly regular vortex lattices have  been heavily
studied both experimentally and theoretically \cite{Abo2001, Madison2001,  Lieb1, Fetter2009}. 
On the other hand,  multi-component BECs admit numerous interesting phenomena absent from 
single-component condensates, for example, domain walls, vortons, square vortex lattices and so on; see 
\cite{Hall1998, Ho1996, Chui2000, Jezek2001, Kasamatsu2004, Kasamatsu2005, 
Kobyakov2011} and references 
therein. As the simplest cases, two-component BECs provide a good opportunity to 
investigate the properties of multi-component condensation.

The first experiment of  two-component BECs was carried out in  
$|F=2,m_f=2\rangle$ and $|F=1, m_f=-1\rangle$ hyperfine  states  of $^{87}{\rm Rb}$ \cite{MBGCW}.  
At temperatures $T$ much smaller than the critical temperature $T_c$,  a rotating two-component
BEC  with an external driving field (or an internal Josephson junction)  can be well described by two 
{self-consistent} nonlinear Schr\"{o}dinger equations (NLSEs),  also known as the  coupled  
Gross--Pitaevskii equations (CGPEs).  The dimensionless CGPEs has the following form  
\cite{Kasamatsu2003,  Kasamatsu2005,  Wang2007, Zhang2007, Bao2008, Bao2011, Kobyakov2011}:
\bea\label{CGPEs1}
&&i\fl{\p\psi_1(\bx,t)}{\p t} = \left[-\fl{1}{2}\nabla^2 + V_1(\bx) + (\bt_{11}|\psi_1|^2+\bt_{12}|\psi_2|^2)
-\Og L_z\right]\psi_1 - \lambda \psi_2, \\
\label{CGPEs2}
&&i\fl{\p\psi_2(\bx, t)}{\p t} = \left[-\fl{1}{2}\nabla^2 + V_2(\bx) + (\bt_{21}|\psi_1|^2+\bt_{22}|\psi_2|^2)
-\Og L_z\right]\psi_2 - \lambda \psi_1, \quad \bx\in{\Bbb  R}^d, \ \, t >0.\qquad\quad
\eea
Here, $\bx\in{\mathbb R}^d$  ($d =  2$ or $3$) is the Cartesian coordinate vector, $t$ is the time and
$\psi_j(\bx,t)$ is the complex-valued macroscopic wave function of the $j$th  ($j = 1, 2$) component.
The interaction constants  $\bt_{jk} = \bt_{kj} = {4\pi N a_{jk}}/{a_0}$ (for $j, k  = 1, 2$), where $N$ is
the  total number of atoms in two-component BECs, $a_0$ is the dimensionless spatial unit
and   $a_{jk} = a_{kj}$ represents the $s$-wave scattering lengths between the $j$th and $k$th
components (positive for repulsive interaction and negative  for attractive interaction).
The constant $\lambda$ describes the effective Rabi frequency to realize the internal atomic 
Josephson junction
by a { Raman transition}, $\Og$ represents the  speed of angular momentum rotation and
$L_z = -i(x\p_y - y\p_x)$ is the $z$-component of the angular momentum operator.  The real-valued function
 $V_j(\bx)$ ($j = 1, 2$)  represents the external trapping potential imposed on the $j$th component.  In
 most BEC experiments,  a harmonic potential is used, i.e.,
 \bea\label{potential}
 V_j(\bx) = \fl{1}{2}\left\{\begin{array}{ll}
 \gm_{x,j}^2x^2 + \gm_{y,j}^2y^2, \qquad & d = 2, \\
 \gm_{x,j}^2x^2 + \gm_{y,j}^2y^2 + \gm_{z,j}^2 z^2, \quad \  & d = 3,
 \end{array}\right. \quad \ \  j = 1, 2.
 \eea
The initial conditions of (\ref{CGPEs1})--(\ref{CGPEs2}) are given by
\bea\label{initial}
\psi_j(\bx, 0) = \psi_j^0(\bx), \quad \ \bx\in{\mathbb R}^d, \quad  \ \ j = 1, 2.
\eea

There are two important invariants associated with the CGPEs  in (\ref{CGPEs1})--(\ref{CGPEs2}):  
the {\it total  mass (or normalization)}, i.e.,
\bea\label{Cnorm}
N(t) : = \|\Psi(\cdot, t)\|^2 = N_1(t) + N_2(t)  \equiv \|\Psi(\cdot, 0)\|^2 = 1, \quad \  t\geq 0,
\eea
where $\Psi(\bx,t) = (\psi_1(\bx,t), \psi_2(\bx,t))^T$  and $N_j$(t) is   the mass of the $j$th component at time 
$t \geq 0$,  which is defined by
\bea
N_j(t) := \|\psi_j(\cdot, t)\|^2  = \int_{{\mathbb R}^d}|\psi_j(\bx,t)|^2 d\bx, \quad \ 
t\geq 0, \quad \ \  j = 1, 2,
\eea
and the {\it energy}
\bea\label{Energy0}
E(\Psi(\cdot, t)) &=& \int_{{\mathbb R}^d}\bigg[\sum_{j=1}^2\bigg(\fl{1}{2}|\nabla\psi_j|^2 +V_j(\bx)|\psi_j|^2
+\fl{\bt_{jj}}{2}|\psi_j|^4 -\Og {\rm Re}\left(\psi_j^*L_z\psi_j\right)\bigg)\qquad\nn\\
&&\qquad +\bt_{12}|\psi_1|^2|\psi_2|^2- 2\lambda{\rm Re}(\psi_1\psi_2^*)\bigg]d\bx  = E(\Psi(\cdot, 0)), \quad \ 
t\geq 0,
\eea
where ${\rm Re}(f)$ and $f^*$ represent the real part and the conjugate of a function $f$, respectively.
In fact, if there is no external driving filed (i.e., $\lambda = 0$ in (\ref{CGPEs1})--(\ref{CGPEs2})),  the
{\it mass of each component} is also conserved, i.e.,  $N_j(t) =  N_j(0)$ ($j = 1, 2$) for $t\geq0$. 
These invariants can be used, in particular,  as benchmarks and validation of numerical algorithms for solving 
the CGPEs  (\ref{CGPEs1})--(\ref{CGPEs2}).

Many numerical methods have been proposed to study the dynamics of the non-rotating two-component 
BECs, i.e.,  when $\Og=0$,  with/without external driving field \cite{Bao2004, Sepulveda2008,JGSGZ,CW}.   
Compared to other methods, the time-splitting pseudo-spectral method in  \cite{Bao2004}
is one of the most successful methods. It has spectral order of accuracy in space and can be easily 
implemented, i.e., they can achieve both the accuracy and efficiency.   However,  the 
appearance of the angular rotational term hinders the direct application of those methods to study the rotating 
two-component BECs.  Recently,  several numerical methods were proposed for simulating the dynamics of 
rotating two-component BECs \cite{Zhang2007, Wang2007,Bao2008, Corro2009,  Hsueh2011, Jin2013}.  
For example,  in \cite{Zhang2007}, a pseudo-spectral type method was proposed by reformulating the 
problem in two-dimensional polar coordinates or three-dimensional cylindrical coordinates. While 
in \cite{Wang2007}, the authors designed a time-splitting alternating direction implicit method, where  
the angular rotation term is treated in $x$- and $y$-directions separately.
Although these methods
have higher spatial accuracy compared to those finite difference/element methods, 
they have their own limitations.  The method in  \cite{Zhang2007}  is only of second-order or fourth-order in the radial
direction, while the implementation of the method in \cite{Wang2007} could become quite involved. 
One possible approach to overcome those limitations is to relax the constrain of the rotational term, 
which is the  main aim of this paper.  In this paper, we propose a simple and efficient numerical method to 
solve the CGPEs (\ref{CGPEs1})--(\ref{CGPEs2}).  
The main merits of our method are:  (i) Using a  rotating Lagrangian coordinate transform, 
we reformulate the  original CGPEs  in  (\ref{CGPEs1})--(\ref{CGPEs2}) to one without angular momentum 
rotation term. Then, the time-splitting pseudo-spectral method designed for the non-rotating BECs,
which are of spectral order accuracy in space and easy to implemented, can be directly applied to solve 
the CGPEs in new coordinates. Moreover,  (ii) our method solves the CGPEs in two splitting steps 
instead of three steps in literature \cite{Zhang2007, Wang2007}, which makes our method more efficient.

The paper is organized as follows. In Section \ref{section2}, we introduce a rotating Lagrangian
coordinate and then cast the CGPEs (\ref{CGPEs1})--(\ref{initial}) in the new
coordinate system.  A simple and efficient numerical method  is introduced
to discretize the CGPEs under a rotating Lagrangian coordinate in Section  \ref{section3}, which is
subsequently generalized in Section \ref{section4} to solve the VGPEs for multi-component BECs. 
To test its performance, we compare our method with those reported in literature and apply it to 
study the dynamics  of rotating two-component BECs  in Section \ref{section5}. In  Section 
\ref{section6},  we make some concluding remarks.

\section{CGPEs under a rotating Lagrangian coordinate}
\setcounter{equation}{0}
\label{section2}

In this section, we first introduce a rotating Lagrangian coordinate and then reformulate the CGPEs
(\ref{CGPEs1})--(\ref{initial}) in the new coordinate system.  In the following, we will always 
refer the Cartesian coordinates  $\bx$ as the {\it  Eulerian coordinates}.  For any time $t\geq 0$, let
${\bf A}(t)$ be an orthogonal rotational matrix  defined as \cite{Garcia2001, Antonelli2012, Bao2013}
\bea\label{Amatrix}
{\bf A}(t)=\left(\begin{array}{cc}
\cos(\Omega t) & \sin(\Omega t) \\
-\sin(\Omega t) & \cos(\Omega t)
 \end{array}\right),  \quad\ \  \mbox{if \ \ $d = 2$,} \qquad \quad  \ \
         \eea
         and
         \bea\label{Matrixa}
{\bf A}(t)=\left(\begin{array}{ccc}
         \cos(\Omega t) & \sin(\Omega t) & 0 \\
         -\sin(\Omega t) & \cos(\Omega t) & 0 \\
         0 & 0  & 1
         \end{array}\right), \quad\ \  \mbox{if \ \ $d = 3$.} \qquad
\eea
It is easy to verify that ${\bf A}^{-1}(t) ={\bf A}^T(t)$ for any $t\ge0$ and ${\bf A}(0) = {\bf I}$  
with ${\bf I}$ the identity matrix.  
Referring the Cartesian coordinates $\bx$ as the {\it Eulerian coordinates}, we introduce the
 {\it rotating Lagrangian coordinates} $\tbx$ as
\bea\label{transform}
\tbx={\bf A}^{-1}(t) \bx={\bf A}^T(t)\bx \quad \Leftrightarrow \quad \bx= {\bf A}(t){\tbx},   \quad\  \bx\in {\mathbb R}^d,
\quad t\geq 0.
\eea
and reformulate  the wave functions  $\psi_j(\bx,t)$ in the new coordinates as $\phi_j(\tbx, t)$
\bea\label{transform79}
\phi_j(\tbx, t):=\psi_j(\bx, t)= \psi_j\left({\bf A}(t){\tbx},t\right), \quad \ \tbx\in {\mathbb R}^d, \quad t\geq0, \quad\ 
\ j = 1, 2.
\eea
We see that when $d = 3$,  the transformation  in (\ref{transform}) does not  change  the coordinate  
in $z$-direction, that is,  $\tz = z$ and the coordinate transformation essentially occurs only in $xy$-plane 
for any $t\geq0$. Fig. \ref{rot-axis} 
 illustrates  the geometrical relation between the $xy$-plane in the Eulerian coordinates  and  $\tilde{x}\tilde{y}$-plane  in the rotating
 Lagrangian coordinates for $\Omega>0$.
{\color{blue}\begin{figure}[h!]
\centerline{
\includegraphics[height=5.996cm,width=7.96cm]{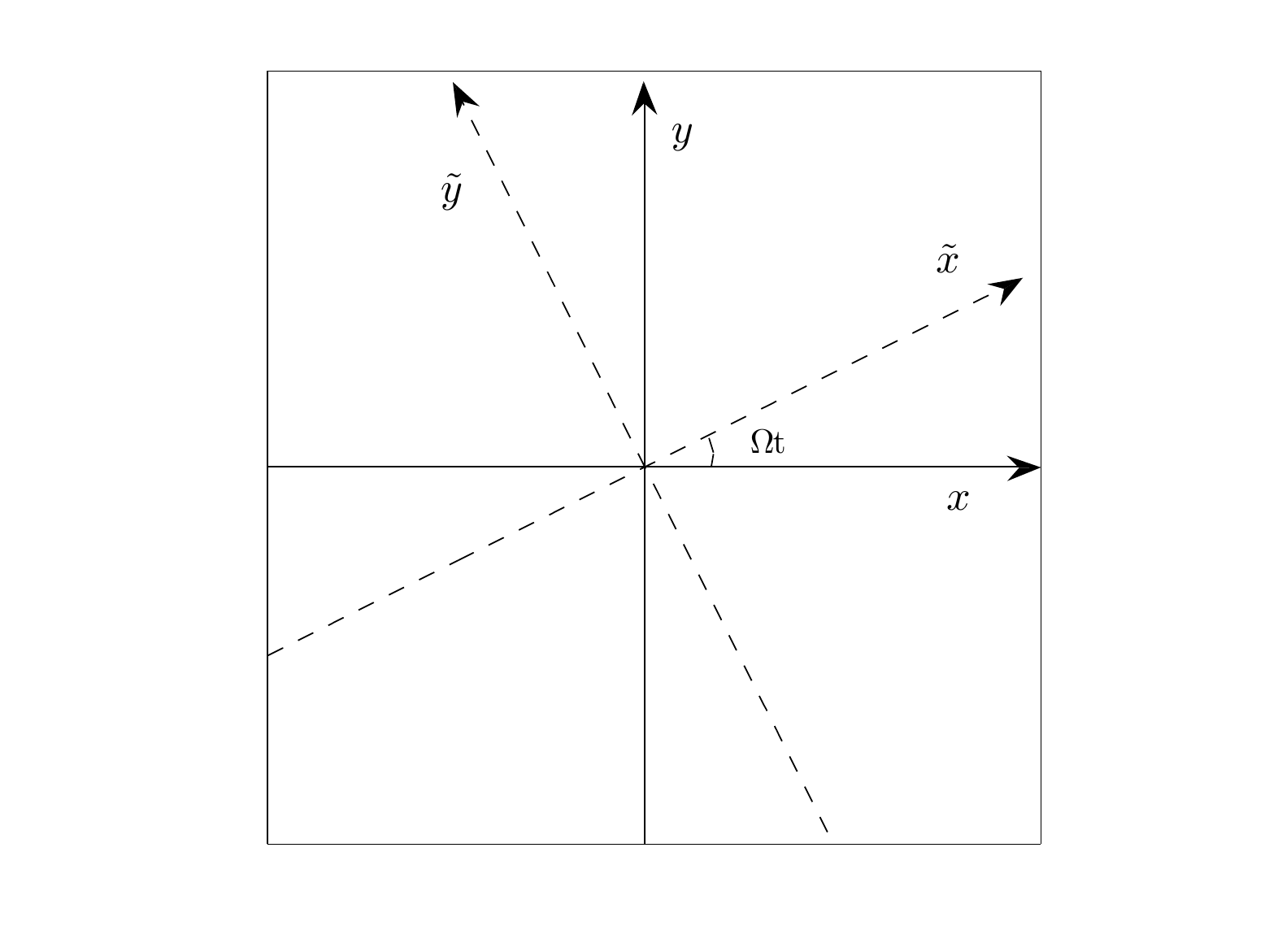}}
\caption{Eulerian (or Cartesian) coordinates $(x,y)$ (solid) and rotating Lagrangian
coordinates $(\tilde{x}, \tilde{y})$ (dashed) in two dimensions for $\Og > 0$ at a fixed $t\ge0$. }\label{rot-axis}
\end{figure}}
Furthermore,  it is easy to see that when $t = 0$ the rotating Lagrangian coordinates $\tbx$ 
become exactly the same as the  Eulerian coordinates $\bx$,  i.e., $\tbx\equiv \bx$ when 
$t = 0$. Note that we assume that  $\Og\neq 0$ in this paper.

From (\ref{transform})--(\ref{transform79}), we  obtain that
\beas\label{operator1}
&&\p_t\phi_j(\tbx,t) =\p_t\psi_j(\bx, t)  + \nabla\psi_j(\bx,t)\cdot\left(\dot{\bf A}(t) \tbx\right)
= \p_t\psi_j(\bx,t)- \Og(x\p_y-y\p_x)\psi_j(\bx,t),\qquad\quad\\
\label{operator2}
&&\nabla\phi_j(\tbx,t) = {{\bf A}^{-1}}(t)\nabla\psi_j(\bx,t) , \quad \nabla^2\phi_j(\tbx, t)
= \nabla^2\psi_j(\bx,t),\quad \ \bx\in{\mathbb R}^d, \quad t\geq0, \quad \  j= 1, 2.
\eeas
 Substituting the above derivatives  into (\ref{CGPEs1})--(\ref{CGPEs2}) gives  the  following
$d$-dimensional CGPEs in the rotating Lagrangian coordinates $\tbx$:
\bea\label{NGPEs1}
&&i\fl{\p\phi_1(\widetilde{\bx}, t)}{\p t} = \left[-\fl{1}{2}\nabla^2 + W_1(\widetilde{\bx}, t) +
(\bt_{11}|\phi_1|^2 + \bt_{12}|\phi_2|^2) \right]\phi_1-
\lambda\phi_2, \qquad\qquad \\
\label{NGPEs2}
&&i\fl{\p \phi_2(\widetilde{\bx}, t)}{\p t} = \left[-\fl{1}{2}\nabla^2 + W_2(\widetilde{\bx}, t) +
(\bt_{21}|\phi_1|^2 + \bt_{22}|\phi_2|^2) \right]\phi_2 -
 \lambda\phi_1, \quad \ \tbx\in{\mathbb R}^d, \quad t >0. \qquad\qquad
\eea
The corresponding initial conditions are
\be
\label{Iinitial}
\phi_j(\tbx, 0) := \phi_j^0(\tbx) =  \psi_j(\bx,0)=\psi_{j}^0(\bx), \quad \ \tbx=\bx\in{\Bbb R}^d.
\ee
In (\ref{NGPEs1})--(\ref{NGPEs2}),  $W_j(\tbx,t)$ ($j = 1, 2$) denotes the effective potential of the
$j$th component, which  is obtained from
\bea\label{wtbx534}
W_j(\tbx,t) = V_j(A(t)\tbx), \quad  \ \tbx\in {\mathbb R}^d, \quad t\geq0, \quad \ \ j = 1, 2.
\eea
In particular, if $V_j(\bx)$ is a harmonic  potential as defined in (\ref{potential}),
then $W_j(\tbx, t)$ has the form
\beas\label{wpotential}
&&W_j(\tbx, t) \\
&&\quad = \fl{\gm_{x,j}^2+\gm_{y,j}^2}{4}
(\tx^2 + \ty^2) + \fl{\gm_{x,j}^2 -\gm_{y,j}^2}{4}\left[(\tx^2-\ty^2)\cos(2\Og t) + 2\tx\ty\sin(2\Og t)\right] +
\left\{\begin{array}{ll} 0
& d= 2, \\
\fl{1}{2}\gm_{z,j}^2\tz^2, & d= 3, \\
\end{array}\right.\qquad 
\eeas
for $j = 1, 2$. Hence, when the external harmonic potentials are radially symmetric in two 
dimensions  (2D)  or
cylindrically symmetric in three dimensions (3D), i.e., $\gm_{x, j}  = \gm_{y, j} :=\gm_{r,j}$, 
 the potential
\bea\label{poten11}
W_j(\tbx, t) = V_j(\tbx), \quad\  \tbx\in{\mathbb R}^d, \quad t\geq0, \quad \  j = 1, 2,
\eea
become  time-independent.

In rotating Lagrangian coordinates, the wave functions $\Phi(\tbx, t) = (\phi_1(\tbx,t), \phi_2(\tbx,t))^T$
satisfy the normalization
\bea\label{Nnorm}
{\widetilde N}(t) : = \|\Phi(\cdot, t)\|^2 = {\widetilde N}_1(t) + {\widetilde N}_2(t)  \equiv \|\Phi(\cdot, 0)\|^2 = 1, 
\quad\  t\geq 0,
\eea
where  ${\widetilde{N}}_j(t)$ is the mass of the $j$th component 
at time $t \geq 0$, i.e.,
\bea
\widetilde{N}_j(t) := \|\phi_j(\cdot, t)\|^2 = \int_{{\mathbb R}^d}|\phi_j(\tbx,t)|^2 d\tbx = N_j(t), \quad
t\geq 0, \quad  \ \ j = 1, 2.
\eea
Similarly, when $\lambda = 0$ the mass of each component is also conserved, i.e., 
$\widetilde{N}_j(t) \equiv \widetilde{N}_j(0)$ for $t \geq0$ and $j = 1, 2$.
The energy associated with the CGPEs (\ref{NGPEs1})--(\ref{NGPEs2}) is 
\bea\label{Energy}
&&\widetilde{E}(\Phi(\cdot, t))=\int_{{\mathbb R}^d}\bigg[\sum_{j=1}^2\bigg(\fl{1}{2}|\nabla\phi_j|^2
+W_j(\tbx, t)|\phi_j|^2 - \int_0^t |\phi_j|^2\p_\tau W_j(\tbx, \tau)d\tau + \fl{\bt_{jj}}{2}|\phi_j|^4\bigg) \nn\\
&&\qquad \qquad \qquad \quad \  \ +\bt_{12}|\phi_1|^2|\phi_2|^2 - 2\lambda{\rm Re}(\phi_1\phi_2^*)
\bigg]
d\tbx  = \widetilde{E}(\Phi(\cdot, 0)), \quad \ t\geq 0.
\eea
As we see in (\ref{poten11}), when $\gm_{x, j}  = \gm_{y, j}$  ($j = 1, 2$),  the potential $W_j(\tbx, t)$ is 
time-independent, which implies that the term of $\int_0^t|\phi_j|^2\p_\tau W_j(\tbx, \tau)d\tau \equiv 0$
in this case.

Compared to (\ref{CGPEs1})--(\ref{CGPEs2}),  the CGPEs (\ref{NGPEs1})--(\ref{NGPEs2}) in
rotating Lagrangian coordinates does not have the angular momentum rotational term, which eliminates 
the difficulties in discretizing the CGPEs and allows us to develop an efficient spectral method to 
solve (\ref{NGPEs1})--(\ref{NGPEs2}).

\section{Numerical method}
\setcounter{equation}{0}
\label{section3}

In this section, we present a time-splitting spectral method to study the dynamics of rotating
two-component BECs.  To the best of our knowledge,  so far all numerical methods computing dynamics  
of rotating two-component BECs in literature are based on discretizing the CGPEs 
(\ref{CGPEs1})--(\ref{CGPEs2})  in Eulerian coordinates {\cite{Zhang2007, Wang2007,  Bao2008}. } 
However, the appearance of rotating angular momentum term in Eulerian coordinates makes it very challenging
to develop an efficient methods with higher accuracy but less computational efforts.  In the following,
instead of simulating (\ref{CGPEs1})--(\ref{CGPEs2}) in Eulerian coordinates,
we solve the CGPEs (\ref{NGPEs1})--(\ref{NGPEs2}) in rotating Lagrangian coordinates. Hence,
we avoid the discretization of the angular rotational term and it makes numerical method simpler 
and more efficient than those reported in {\cite{Wang2007, Zhang2007}.}

In practical computations, we truncate the problem (\ref{NGPEs1})--(\ref{NGPEs2}) into a bounded
computational domain ${\mathcal D}\subset{\Bbb R}^d$ and consider
\bea\label{DCGPEs1}
&&i\p_t \phi_1(\widetilde{\bx}, t) = \left[-\fl{1}{2}\nabla^2 + W_1(\widetilde{\bx}, t) +
(\bt_{11}|\phi_1|^2 + \bt_{12}|\phi_2|^2) \right]\phi_1 - \lambda\phi_2, \qquad\qquad \\
\label{DCGPEs2}
&&i\p_t \phi_2(\widetilde{\bx}, t) = \left[-\fl{1}{2}\nabla^2 + W_2(\widetilde{\bx}, t) +
(\bt_{21}|\phi_1|^2 + \bt_{22}|\phi_2|^2) \right]\phi_2 - \lambda\phi_1,
\quad\  \tbx\in{\mathcal D}, \quad t>0, \qquad\qquad
\eea
along with the initial conditions
\bea\label{Dinitial}
\phi_j(\tbx, t) = \phi_j^0(\tbx), \quad \  \tbx\in\overline{\mathcal D}, \qquad \mbox{with} \quad \int_{\overline{\mathcal D}}
\left(|\phi_1^0(\tbx)|^2 + |\phi_2^0(\tbx)|^2\right)d\tbx = 1.
\eea
The following homogeneous Dirichlet boundary conditions are considered here,  i.e.,
\bea\label{DBC}
\phi_j(\tbx, t) = 0, \quad\  \tbx\in\p{\mathcal D},\quad t>0, \quad \ \ j = 1, 2.
\eea
Due to the confinement of the external potential $W_j(\bx,t)$ and conservation of the normalization 
(\ref{Nnorm}) and energy (\ref{Energy}), the wave function $\phi_j(\tbx,t)$ vanishes as 
$|\tbx|\rightarrow\infty$.  Hence, it is natural to impose homogeneous Dirichlet boundary conditions to 
the truncated problem  (\ref{DCGPEs1})--(\ref{Dinitial}). The use of more sophisticated 
boundary conditions for more generalized cases, e.g., 
absence of trapping potential, is an interesting
topic that remains to be examined in the future \cite{ABK, BT}.
In practical simulations, 
the computational domain ${\mathcal D}\subset{\mathbb R}^d$ is chosen as ${\mathcal D} = 
[a, b]\times[c, e]$
if $d = 2$ and ${\mathcal D} = [a, b]\times[c, e] \times[f, g]$ if $d = 3$. Moreover, 
we use sufficiently large domain  ${\mathcal D}$ to ensure the homogeneous
Dirichlet boundary conditions do not introduce aliasing error. Usually, the diameter  of the bounded
computational domain depends on the problem. In general, it should be larger than the ``Thomas-Fermi
radius" \cite{Bao2005, Bao2006}. 

\subsection{Time-splitting method}
\label{section3-1}

In the following, we use the time-splitting method to discretize the problem (\ref{DCGPEs1})--(\ref{DBC})
in time.  To do it,  we choose a time step $\Dt  t> 0$ and define time sequence $t_n = n\Dt t$ for $n = 0, 1,
\ldots$.   Then from time $t = t_n$  to $t = t_{n+1}$, we numerically solve the CGPEs
(\ref{DCGPEs1})--(\ref{DCGPEs2}) in two steps, i.e., solving
\bea\label{step1}
i\fl{\p\phi_j(\tbx, t)}{\p t} = -\fl{1}{2}\nabla^2\phi_j(\tbx, t) -  \lambda\phi_{(3-j)}(\tbx, t), \quad\  \tbx\in{\mathcal D},\quad 
t_n\leq t\leq t_{n+1}, \quad \ \ j = 1, 2,
\eea
and
\bea\label{step2}
i\fl{\p\phi_j(\tbx, t)}{\p t} =  \Big(W_j(\tbx,t) + \sum_{k = 1}^2 \bt_{jk}|\phi_k|^2\Big)\phi_j(\tbx, t),  \quad\
\tbx\in{\mathcal D},\quad t_n\leq t\leq t_{n+1}, \quad\ \ j = 1, 2.
\eea
In fact, Eq. (\ref{step1}) is  coupled linear Schr\"{o}dinger equations and its discretization will be discussed 
later.

We notice that  in (\ref{step2}), both $|\phi_1(\tbx, t)|$ and $|\phi_2(\tbx, t)|$ are invariants
in time $t$, i.e., $|\phi_j(\tbx, t)| = |\phi_j(\tbx, t_n)|$ ($j = 1, 2$) for any $t\in[t_n,  t_{n+1}]$. Thus, for time 
$t\in[t_n,  t_{n+1}]$, (\ref{step2})  is equivalent to
\bea\label{step20}
i\p_t\phi_j(\tbx, t) =  \Big(W_j(\tbx,t) + \sum_{k = 1}^2 \bt_{jk}|\phi_k(\tbx, t_n)|^2\Big)\phi_j(\tbx, t),
\quad\  \tbx\in{\mathcal D},\quad \ \ j = 1, 2.
\eea
Integrating (\ref{step20}) exactly in time leads to the solution of (\ref{step2}), i.e.,
\bea\label{step2sol}
\phi_j(\tbx, t) = \phi_j(\tbx, t_n)\exp\Big[-i\Big((t-t_n)\sum_{k = 1}^2 \bt_{jk}|\phi_k(\tbx, t_n)|^2
 + \int_{t_n}^t  W_j(\tbx, \tau)d\tau\Big)\Big], \quad\ \  j = 1, 2,
\eea
for $\tbx\in{\mathcal D}$ and $t\in[t_n, t_{n+1}]$.

\vskip 20pt
\begin{remark}
If $V_j(\bx)$  ($j = 1, 2)$ is a harmonic potential as defined in (\ref{potential}), then the integral in
(\ref{step2sol}) can  be evaluated analytically, i.e.,
\bea
\int_{t_n}^t W_j(\tbx, \tau)d\tau = \fl{(\gm_{x,j}^2+\gm_{y,j}^2)(\tx^2+\ty^2)}{4}(t - t_n) + U(\tbx, t) + 
\left\{\begin{array}{ll} 0, &d = 2, \\
\fl{1}{2}\gm_{z,j}^2\tz^2(t-t_n),  \quad &d = 3,
\end{array}\right.
\eea
where
\beas
U(\tbx,t) &=& \fl{(\gm_{x,j}^2-\gm_{y,j}^2)}{4}\int_{t_n}^t \left[(\tx^2-\ty^2)\cos(2\Og\tau)
+2\tx\ty\sin(2\Og\tau)\right]d\tau\\
&=&\fl{(\gm_{x,j}^2-\gm_{y,j}^2)(\tx^2-\ty^2)}{8\Og}[\sin(2\Og t) - \sin(2\Og t_n)] -
\fl{(\gm_{x,j}^2-\gm_{y,j}^2)\tx\ty}{4\Og}[\cos(2\Og t)-\cos(2\Og t_n)].
\eeas
Typically, when $\gm_{x,j} = \gm_{y,j}$ ($j = 1, 2$), we have $U(\tbx, t) \equiv 0$.

For a general potential $V_j(\bx)$, if the integral in (\ref{step2sol}) can not be found analytically,
numerical quadratures such as Trapezoidal rule or Simpson's rule can be used to calculate its approximation
\cite{Bao2006, Bao2013}.
\end{remark}

\subsection{Discretization of coupled linear Schr\"{o}dinger equations}
\label{section3-2}

In the following, we first introduce a linear transformation of the wave functions $\phi_j(\tbx,t)$ 
($j = 1, 2$) such that the coupled linear Schr\"{o}dinger equations become independent of each 
other. Then we describe the  sine pseudospectral discretization in two-dimensional case. Its 
generalization to three dimensions is  straightforward.

Let the  matrix
\bea
{\bf P} = \left(\begin{array}{rr}
1 & 1\\
1 & -1\end{array}\right),
\eea
and denote
\bea\label{trans}
\left(\begin{array}{c}
\varphi_1(\tbx,t)\\
\varphi_2(\tbx,t)
\end{array}\right)= {\bf P} \left(\begin{array}{c}
\phi_1(\tbx, t)\\
\phi_2(\tbx, t)
\end{array}\right) = \left(\begin{array}{c}
\phi_1(\tbx, t) + \phi_2(\tbx, t)\\
\phi_1(\tbx, t) - \phi_2(\tbx, t)
\end{array}\right), \quad \ \tbx\in{\mathbb R}^d, \quad t\geq0.
\eea
Combining (\ref{step1}) and (\ref{trans}),  we obtain the following equations for $\varphi_j(\tbx,t)$
($j = 1, 2$):
\bea\label{step10}
&&i\p_t\varphi_1(\tbx,t) =  -\fl{1}{2}\nabla^2\varphi_1(\tbx,t) - \lambda\varphi_1(\tbx, t),\\
\label{step11}
&&i\p_t\varphi_2(\tbx,t) =  -\fl{1}{2}\nabla^2\varphi_2(\tbx,t) + \lambda\varphi_2(\tbx, t), \quad\  \tbx\in{\mathcal D}, \ \ \
t_n\leq t\leq t_{n+1}.
\eea
It is easy to see that the functions $\varphi_1$ and $\varphi_2$  are independent in (\ref{step10})--(\ref{step11}), which allows us to solve them separately.

Choose two even integers $J, K > 0$ and  denote the index set
\beas
{\mathcal T}_{JK} = \left\{(p, q)\,|\, 1\leq p\leq J-1, \ \ 1\leq q\leq K-1\right\}.
\eeas
Define the function
\beas
U_{pq}(\tbx) = \sin(\mu_p^x(\tx-a))\sin(\mu_q^y(\ty-c)),  \quad \
\tbx=(\tx, \ty)^T \in{\mathcal D}, \quad\ \  (p, q)\in{\mathcal T}_{JK}.
\eeas
Assume that
\bea\label{anz}
\varphi_j(\tbx, t)  =  \sum_{p = 1}^{J-1}\sum_{q = 1}^{K-1} \widehat{\varphi}_{j, pq}(t)\, U_{pq}(\tbx),
\quad\  \tbx\in{\mathcal D},\quad t \in[t_n, t_{n+1}], \quad\ \  j = 1, 2,
\eea
where  $\widehat{\varphi}_{j,pq}(t)$ is the discrete 
sine transform of $\varphi_{j}(\tbx,t)$ corresponding to frequencies $(p,q)$ and
\beas
\mu_p^x = \fl{p\pi}{b-a}, \quad\  \mu_q^y = \fl{q\pi}{e-c}, \qquad (p, q)\in{\mathcal T}_{JK}.
\eeas
Substituting (\ref{anz}) into (\ref{step10})--(\ref{step11})  leads to
\bea\label{sol11}
&&\widehat{\varphi}_{1, pq}(t) = \widehat{\varphi}_{1,pq}(t_n) \exp\Big[-i\Big(\fl{(\mu_p^x)^2
+(\mu_q^y)^2}{2}-\lambda\Big)(t-t_n)\Big],  \\
\label{sol12}
&&\widehat{\varphi}_{2, pq}(t) = \widehat{\varphi}_{2,pq}(t_n) \exp\Big[-i\Big(\fl{(\mu_p^x)^2
+(\mu_q^y)^2}{2}+\lambda\Big)(t-t_n)\Big], \quad (p,q)\in{\mathcal T}_{JK}, \ \ t\in[t_n, t_{n+1}].
\quad\qquad
\eea
Combining (\ref{sol11})--(\ref{sol12}) and  (\ref{anz}) and noticing the linear transformation in
  (\ref{trans}), we obtain a sine pseudospectral approximation  to (\ref{step1}), i.e.,
\bea\label{sol-step0}
&&\phi_j(\tbx,t)=\sum_{p = 1}^{J-1}\sum_{q = 1}^{K-1}\left[\cos(\lambda(t-t_n))
\widehat{\phi}_{j,pq}(t_n)+i\sin(\lambda(t-t_n))\widehat{\phi}_{(3-j),pq}(t_n)\right]
 \eta_{pq}(t)\,U_{pq}(\tbx),\qquad\quad
\eea
for $j = 1, 2$, where
\beas
\eta_{pq}(t) = \exp\Big[-\fl{i}{2}\left((\mu^x_p)^2+(\mu^y_q)^2\right)(t-t_n)\Big], \qquad
(p, q) \in{\mathcal T}_{JK}.
\eeas
In (\ref{sol-step0}), $\widehat{\phi}_{j,pq}(t)$ is the discrete sine transform of $\phi_j(\tbx, t)$ 
($j = 1, 2$) corresponding to the frequency $(p, q)$. We remark here that although the solution 
(\ref{sol-step0}) is found via (\ref{sol11})--(\ref{sol12}), (\ref{anz}) and (\ref{trans}), in practice 
we only need to compute  $\widehat{\phi}_{1,pq}(t)$ and $\widehat{\phi}_{2,pq}(t)$ to obtain 
(\ref{sol-step0}).

\subsection{Implementation of the method}
\label{section3-3}

For convenience of the readers,  in the following we will summarize our method and describe its 
implementation.  For simplicity of notations, the method will only be presented in two-dimensional case.
Choose spatial mesh sizes $h_{\tx} = (b-a)/J$ and $h_{\ty} = (e-c)/K$ in $\tx$- and $\ty$-directions, respectively.  
Define
\beas
\tx_s = a+sh_{\tx},  \quad 0\leq s\leq J;\qquad \ty_l = c+lh_{\ty},\quad 0\leq l\leq K.
\eeas
Let $\phi_{j, sl}^n$ denote the numerical approximation to $\phi_{j}(\tx_k, \ty_l, t_n)$.
From $t = t_{n}$ to $t = t_{n+1}$, we use the second-order Strang splitting method 
\cite{Strang1968, Glowinski1989, Bao2005} to combine
the two steps in (\ref{step1})  and (\ref{step2}), i.e,
\bea\label{ss1}
&&\phi_{j,sl}^{(1)} = \phi_{j,sl}^n \exp\Big[-i\Big(\fl{\Dt t}{2}\sum_{k=1}^2\bt_{jk}|\phi_{k,sl}^n|^2 
+ \int_{t_n}^{t+\Dt t/2}  W_j(\tx_s, \ty_l, \tau)d\tau\Big)\Big],\\
\label{ss3}
 &&\phi_{j,sl}^{(2)}=   \sum_{p = 1}^{J-1}\sum_{q = 1}^{K-1} e^{-i\fl{\Dt t}{2}[(\mu_p^x)^2+(\mu_q^y)^2]}
 \Big[\cos(\lambda\Dt t)\widehat{\phi}_{j,pq}^{(1)}
 +i\sin(\lambda\Dt t)\widehat{\phi}_{(3-j),pq}^{(1)}\Big]
 \sin\left(\fl{sp\pi}{J}\right)\sin\left(\fl{lq\pi}{K}\right), \qquad \\
 \label{ss2}
 &&\phi_{j,sl}^{n+1} = \phi_{j,sl}^{(2)} \exp\Big[-i\Big(\fl{\Dt t}{2}\sum_{k=1}^2\bt_{jk}|\phi_{k,sl}^{(2)}|^2 
 + \int_{t+\Dt t/2}^{t_{n+1}} W_j(\tx_s, \ty_l, \tau)d\tau\Big)\Big],\qquad j = 1, 2,
\eea
for $0\leq s\leq J$,  $0\leq l\leq K$ and $n = 0, 1, \ldots$.  At $t = 0$, the initial conditions (\ref{Dinitial}) are 
discretized as
\bea\label{ss0}
\phi_{j,sl}^0 = \phi_j^0(\tx_s, \ty_l), \quad  \ 0\leq s\leq J, \quad 0\leq l\leq K, \qquad j = 1, 2,
\eea

Our method described in (\ref{ss1})--(\ref{ss0}) is explicit and it is easy to implement. 
Furthermore,  the memory cost is $O(JK)$ and the computational cost per time step is 
$O(JK\ln (JK))$ if a 2D CGPEs is solved.   In 3D case, the memory cost and the computational 
cost per time step are $O(JKL)$ and $O(JKL\ln(JKL))$, respectively, where the even integer 
$L > 0$ and $L+1$ is the number of grid points in $z$-direction in 3D.

\begin{remark}\label{Remark2}
The solutions $\phi_{1,sl}^{n+1}$ and $\phi_{2,sl}^{n+1}$ obtained from (\ref{ss1})--(\ref{ss2})
are  grid functions on the bounded computational domain ${\mathcal D}$ in  rotating
Lagrangian coordinates.  To obtain the wave functions $\psi_1(\bx, t_{n})$ and $\psi_2(\bx,t_n)$ 
satisfying the CGPEs  (\ref{CGPEs1})--(\ref{CGPEs2})  over a set of fixed grid points in the Eulerian 
coordinates $\bx$, we can  use the standard Fourier/sine interpolation operators from the discrete 
numerical solution $\phi(\tbx,t_n)$ to construct an interpolation continuous function over 
${\mathcal D}$ \cite{Boyd1992, Shen}.
\end{remark}

\begin{remark}\label{Remark3}
If the potential $V_j(\bx)$ in (\ref{potential}) is replaced by a time-dependent potential, e.g.,
$V_j(\bx,t)$, the rotating Lagrangian coordinate transformation and the numerical method
are still valid provided that we replace $W_j(\tbx,t)$ in (\ref{wtbx534}) by
$W_j(\tbx,t) = V_j(A(t)\tbx,t)$ for $\tbx\in {\mathbb R}^d$ and $t\ge0$.
\end{remark}

\section{Extension to rotating multi-component BECs}
\setcounter{equation}{0}
\label{section4}

In Sections \ref{section2}--\ref{section3}, we presented an efficient and accurate numerical method 
to compute the dynamics of rotating two-component BECs with internal Josephson junction. 
In fact, this method can be easily generalized to solve the vector Gross--Pitaevskii equations (VGPEs)
with an angular momentum rotation term and an external driving field, which describes the 
dynamics of  rotating multi-component BECs \cite{Bao2004, CLL, Liu, Lieb1}.
 
Suppose there are $M \geq 2$ species in  multi-component BECs. Denote 
the complex-valued macroscopic wave function for the $j$th component 
as $\psi_j(\bx,t)$ for $j=1,\ldots,M$.  
Let $\Psi(\bx,t)=\left(\psi_1(\bx,t),\ldots,\psi_M(\bx,t)\right)^T$.  Then the evolution of the wave function 
$\Psi(\bx,t)$  is governed by the following self-consistent VGPEs 
\cite{Bao2004, CLL, Liu, Lieb1, LW}:
  
\begin{equation}\label{VGPEs}
i\fl{\p\Psi(\bx, t)}{\p t} = \left[-\fl{1}{2}\nabla^2 + {\bf V}(\bx) + {\bf F}(\Psi)  
-\Og L_z + g(t){\bf B} \right]\Psi, \quad\  \bx\in{\Bbb  R}^d, \quad t >0
\end{equation}
with the initial conditions
\begin{equation}\label{VGPEs-ini}
\Psi(\bx,0)=\Psi^0(\bx) = (\psi_1^0(\bx), \ldots, \psi_M^0(\bx))^T, \quad\ \  \bx\in {\Bbb  R}^d. 
\end{equation}
The matrix   ${\bf V}(\bx)={\rm diag}\left( V_1(\bx), \, \ldots,\, V_M(\bx) \right)^T$ represents 
the external traping potentials and 
${\bf F}(\Psi)={\rm diag}\left( F_1(\Psi),\, \ldots,\, F_M(\Psi) \right)^T$ with
$$F_j(\Psi)=\sum_{k=1}^{M}\beta_{jk}|\psi_k(\bx,t)|^2,\quad\ \  j=1, \ldots, M,$$ 
where the constant $\beta_{jk}$ describes the interaction strength between the 
$j$th and  $k$th components.  $g(t)$ is a real-valued scalar function and  ${\bf B}$ is 
a real-valued diagonalizable constant matrix.

To solve  (\ref{VGPEs})--(\ref{VGPEs-ini}), similarly we introduce the rotating Lagrangian 
coordinates as defined in (\ref{transform})--(\ref{transform79}) and cast the VGPEs in the new 
coordinates.  Then we truncate it into a bounded computational domain ${\mathcal D}\subset{\mathbb R}^d$ 
and consider  the following VGPEs with homogenous Dirichlet boundary conditions: 
\bea
\label{rot-VGPEs1}
&&i\fl{\p\Phi(\tbx, t)}{\p t} = \left[-\fl{1}{2}\nabla^2 + {\bf W}(\tbx,t) + 
{\bf F}(\Phi) + g(t){\bf B} \right]\Phi, \quad \  \tbx\in{\mathcal D}, \quad t >0, \qquad\qquad\qquad\quad\\
&&
\label{rot-VGPEs2}
 \Phi(\tbx,0)=\Phi^0(\tbx),\quad \  \tbx\in\overline{\mathcal D} \qquad\mbox{and}\qquad \Phi(\tbx,t)=0, 
 \quad \  \tbx\in{\p \mathcal D}, \quad t > 0,
\eea
where   $\Phi(\tbx,t)=\left(\phi_1(\tbx,t),\,\ldots,\,\phi_M(\tbx,t)\right)^T$
and   ${\bf W}(\tbx,t)={\rm diag} \left(W_1(\tbx,t),\, \ldots,\, W_M(\tbx,t) \right)^T$  
with $W_j(\tbx,t) = V_j(A(t)\tbx)$ for $j = 1, \ldots, M$. From time $t = t_n$ to $t = t_{n+1}$,  we split the 
VGPEs (\ref{rot-VGPEs1}) into two subproblems and solve
\begin{equation}\label{SP-rotVGPEs1}
i\p_t\Phi(\tbx, t) = \left[-\fl{1}{2}\nabla^2 +  g(t) {\bf B} \right]\Phi, \quad \  \tbx\in{\mathcal D}, \quad t_n \leq 
t \leq t_{n+1}
\end{equation}
for a time step of length $\Dt t$, followed by solving 
\begin{equation}\label{SP-rotVGPEs2}
i\p_t\Phi(\tbx, t) = \left[ {\bf W}(\tbx,t) + {\bf F}(\Phi)  \right]\Phi,\quad \  \tbx\in{\mathcal D}, \quad t_n \leq 
t \leq t_{n+1}
\end{equation}
for the same time step.

Equation (\ref{SP-rotVGPEs2}) can be integrated exactly in time and the solution is 
\bea\label{s1}
\phi_j(\tbx, t) = \phi_j(\tbx,t_n) \exp\Big[-i\Big((t-t_n)F_j(\Phi(\tbx, t_n))+\int_{t_n}^tW_j(\tbx,\tau)
d\tau\Big)\Big],\quad \ 
j = 1, \ldots, M
\eea
for $\tbx\in{\mathcal D}$ and $t\in[t_n, t_{n+1}]$. On the other hand,  since ${\bf B}$ is a  diagonalizable 
matrix, there exists a matrix  ${\bf D}$ and a diagonal 
matrix ${\bf \Lambda}={\rm diag}\left(\lambda_1, \,\ldots,
\,\lambda_M \right) $  such that ${\bf B}={\bf D}^{-1}{\bf \Lambda}{\bf D}.$   Denote 
\bea\label{Upsilon}\Upsilon(\tbx,t):={\bf D}\Phi(\tbx,t) =\left(\varphi_1(\tbx,t),\, \ldots,\,\varphi_M(\tbx,t)
\right)^T.
\eea
Then from (\ref{SP-rotVGPEs2}),  we obtain 
\bea\label{Upsilon1}
i\p_t \Upsilon(\tbx, t) = \left[-\fl{1}{2}\nabla^2 + g(t){\bf \Lambda}\right]\Upsilon,  \quad \  
\tbx\in{\mathcal D}, \quad t_n \leq 
t \leq t_{n+1}.
\eea
Again we will only presented its solution in 2D case and the generalization to 3D is straightforward.  
Following  the similar procedures in Sec. \ref{section3-2}, i.e., solving (\ref{Upsilon1}) and noticing 
$\Phi = {\bf D}^{-1}\Upsilon$,  we obtain the solution of (\ref{SP-rotVGPEs1}) as
\bea\label{s2}
\Phi(\tbx, t) = \left({\bf D}^{-1}e^{-i{\bf \Lambda}\int_{t_n}^t g(\tau)d\tau}\right)\sum_{p=1}^{J-1}
\sum_{q=1}^{K-1}e^{-\fl{i}{2}\left[(\mu_p^x)^2
+(\mu_q^y)^2\right] (t-t_n)}\left({\bf D}\widehat{\Phi}_{pq}(t_n)\right)U_{pq}(\tbx), \eea
for $\tbx\in{\mathcal D}$ and $t\in[t_n, t_{n+1}]$, where $\widehat{\Phi}_{pq}(t_n) = (
\widehat{\phi}_{1,pq}(t_n), \widehat{\phi}_{2,pq}(t_n),  \ldots, \widehat{\phi}_{M,pq}(t_n))^T$
with $\widehat{\phi}_{j,pq}(t_n)$  ($j = 1, 2, \ldots, M)$ the discrete sine transform of $\phi_j(\tbx,t_n)$
corresponding to the frequency $(p,q)$.
Similarly, we can use the second-order Strange splitting method to combine the above two 
steps and it  can be easily implemented by replacing (\ref{ss1}) and (\ref{ss2}) by (\ref{s1}) and 
 (\ref{ss3}) by  (\ref{s2}). 

\section{Numerical results}
\setcounter{equation}{0}
\label{section5}

In this section, we first test the accuracy of our method presented in Sec. \ref{section3} and compare
its accuracy and efficiency with the method reported in \cite{Wang2007}. Then we apply our method to study the dynamics of vortex lattices and  other properties of rotating two-component BECs.

\subsection{Comparison of methods}

In this section, we test the accuracy and efficiency of our method and compare it  with the time-splitting
alternating direction implicit (TSADI) method proposed in \cite{Wang2007}.  The TSADI method
has spectral accuracy in all spatial directions and thus is more accurate than other methods in 
literature.  
Hence, in the following we only compare our method with the TSADI method
in \cite{Wang2007}.

We solve the two-dimensional (i.e., $d = 2$)  CGPEs  with the  parameters $\Og = 0.4$,
$\lambda  = 1.0$, $\gm_{x,j} = \gm_{y,j} = 1$ (for $j = 1, 2$),  and
\beas
\left(\begin{array}{c c}
 \bt_{11} & \bt_{12} \\
 \bt_{21} & \bt_{22}\end{array}\right) = 50\left(\begin{array}{cc}
 1.03 & 1.0 \\
 1.0 & 0.97 \end{array}\right).
\eeas
The initial conditions are chosen as
\bea
\psi_1^0(\bx) =  \fl{1}{\sqrt{2\pi}}e^{-\fl{x^2+y^2}{2}}, \qquad
\psi_2^0(\bx) =  \fl{1.5^{1/4}}{\sqrt{2\pi}}e^{-\fl{x^2+1.5y^2}{2}}, \qquad \bx\in{\mathbb R^2}.
\eea
We remark here that the TSADI method in \cite{Wang2007}  is different from our method  mainly in
three aspects:  (i) The TSADI method solves the CGPEs  (\ref{CGPEs1})--(\ref{CGPEs2}) in Eulerian 
coordinates.  While our method solves the CGPEs (\ref{NGPEs1})--(\ref{NGPEs2}) in rotating 
Lagrangian coordinates.  (ii) To decouple the nonlinearity and 
internal Josephson junction terms, the TSADI method splits the CGPEs into three steps, while our 
method can solve  the problem in two steps.   (iii) In \cite{Wang2007}, the angular momentum 
rotational term
$-\Og L_z\psi$ is ``split" into two parts using ADI method. In contrast,  in our method we 
use coordinate transformation to eliminate this term and thus avoid to discretize it. For details of the 
TSADI method, we refer readers to \cite{Bao2007, Wang2007}\footnote{See Eq. ({3.14}) in \cite{Wang2007} for more
details of the TSADI method. Notice that  in (3.14) $\psi_{1,jk}^{(3)}$, $\psi_{2, jk}^{(3)}$, 
$\psi_{1,jk}^{(5)}$ and $\psi_{2,jk}^{(5)}$ are mistyped. For example, $\psi_{1,jk}^{(3)}$ and
$\psi_{2, jk}^{(3)}$  should be computed as \cite{Zhang2007, Bao2008}
\bea\label{correction1}
&&\psi_{1,jk}^{(3)} = \cos(\lambda\Dt t/2)\psi_{1,jk}^{(2)} + i\sin(\lambda\Dt t/2)\psi_{2,jk}^{(2)},\qquad\qquad\\
\label{correction2}
&&\psi_{2,jk}^{(3)} = \cos(\lambda\Dt t/2)\psi_{2,jk}^{(2)} + i\sin(\lambda\Dt t/2)\psi_{1,jk}^{(2)}.\qquad\qquad
\eea
Similarly, $\psi_{1,jk}^{(5)}$  and $\psi_{2, jk}^{(5)}$ should also be changed correspondingly. 

}.

Since the TSADI and our method solve the problem in different coordinates,  to 
compare them in a fair way we will use the same spatial mesh size and time step.  In simulations, 
we choose sufficiently large computational domain ${\mathcal D} = [-16,16]^2$ for both methods.
Denote $\phi_j^{(h_{\tx}, h_{\ty}, k)}(t)$ as the numerical approximation of $\phi_j(\tbx, t)$, which is 
obtained by using our method with time step $k$ and spatial 
 mesh size $h_{\tx}$ and $h_{\ty}$.  Similarly, 
let $\psi_j^{(h_x, h_y, k)}(t)$ be the numerical solution of $\psi_j(\bx, t)$ from the TSADI method.
Here we take $h_x = h_y = h_{\tx} = h_{\ty} := h$.  With a slight abuse of notation, 
we let  $\phi_j(t)$ (or $\psi_j(t)$) represent
the numerical solution with very fine mesh size $h = 1/64$ and small time step $k = 0.0001$
and assume it to be sufficiently good representation of the exact solution at time $t$.
Tables \ref{T1}--\ref{T2} show the spatial and temporal errors of two methods, where the errors
are computed by
\beas
\|\Phi(t) - \Phi^{(h,h,k)}(t)\|_{l^2} = \sqrt{\sum_{j=1}^2\|\phi_j(t) - \phi_j^{(h,h,k)}(t)\|^2_{l^2}}
\eeas
for our method and  $\|\Psi(t) - \Psi^{(h, h, k)}(t)\|_{l^2}$ for the TSADI method. In addition, we show the
CPU time consumed by each method for time $t\in[0, 2]$. To calculate the spatial errors in Table
\ref{T1}, we always use a very small time step $k = 0.0001$ so that the errors from time discretization
can be neglected compared to those from spatial discretization. On the other hand, in Table
\ref{T2}, we always use $h = 1/64$ which is the same as those used in obtaining
the `exact' solution, so that one can regard the spatial discretization as  `exact' and the only errors
are from time discretization.

\begin{table}[htb!]
\begin{center}
\begin{tabular}{|c|c|r|c|r|}
\hline
 &\multicolumn{2}{c|}{TSADI method} &
 \multicolumn{2}{c|}{Our method}\\
 \hline
Mesh size $h$ & Error  & Computing time (s.) & Error & Computing time (s.) \\
\hline
1     & 0.8562 & 19.50 & 0.9408&13.16 \\
1/2  & 0.1202 & 81.19 & 0.1202& 53.09\\
1/4  & 6.9425E-4 & 383.74 & 6.8771E-4&  235.54\\
1/8  & 1.1267E-7 & 1627.15 & 3.8578E-8& 1002.43\\
1/16& $<$ 1.0E-8 & 7450.26 & $<$ 1.0E-8& 4502.41\\
 \hline
\end{tabular}
\caption{Spatial discretization errors at time $t = 2$ and the computing time (i.e., CPU time in
second)  spent  by  each method, where the time step $k = 0.0001$ for both methods. }\label{T1}
\end{center}
\end{table}

\begin{table}[htb!]
\begin{center}
\begin{tabular}{|c|c|r|c|r|}
\hline
 &\multicolumn{2}{c|}{TSADI method} &
 \multicolumn{2}{c|}{Our method}\\
 \hline
Time step $k$ & Error  & Computing time (s.)& Error & Computing time (s.) \\
\hline
1/40   &  1.7511E-2 & 649.14 & 1.0164E-2& 415.21\\
1/80   &  4.3444E-3 & 1277.80 & 2.5310E-3& 817.31 \\
1/160 &  1.0839E-3 & 2587.39& 6.3204E-4 & 1620.82 \\
1/320 &  2.7064E-4  & 5033.16 &1.5785E-4 & 3207.06\\
1/640 &  6.7444E-5 & 9956.51 & 3.9339E-5& 6383.44  \\
 \hline
\end{tabular}
\caption{Temporal discretization errors at time $t = 2$ and the computing time (i.e., CPU time in
second) spent  by  each method, where the mesh size $h = 1/64$ for both methods. }\label{T2}
\end{center}
\end{table}

From Tables \ref{T1}--\ref{T2}, we see that both the TSADI method and our method have
the spectral accuracy in space and  the second-order of accuracy in time.  However, for the same
numerical parameters (i.e., $h$ and $k$), the TSADI method is much slower
than our method. Usually, the computing time spent by the TSADI method is around 1.5
times more than that taken by our method. For example, when $h = 1/64$ and $k = 1/160$, the
computing time by the TSADI method is $2587.39$ and our method only needs $1620.82$.
This is mainly caused  by two factors: (i)  The TSADI method splits the spatial operator into the
operators in $x$- and $y$-directions.  Equivalently, it discretizes the $x$- and $y$-direction separately.
While our method treats all spatial directions simultaneously, which saves the time in doing  
discrete sine transform.  (ii) In \cite{Wang2007}, the CGPEs is 
solved by  three splitting steps, i.e., there is  an extra step of $i\p_t\psi_j = -\lambda \psi_{(3-j)}$ 
for $j = 1, 2$ to solve at each time step.  However, in our method we notice that the term of 
$-\lambda\psi_{(3-j)}$ can be combined with the $\Dt\phi_j$ by a linear transformation of the wave 
functions, which avoids introducing the extra step to treat the internal Josephson junction terms.
Hence, our method is more efficient especially in higher dimensions or when more components are 
involved.

The computing time of both methods increases when smaller time step or spatial mesh size are 
used. Especially,  for a fixed time step $k$, if the mesh size $h$ decreases by a factor $\ap_h$, then
the time spent by both methods increases by a factor of $\ap_h^2$.  While for  fixed mesh size $h$,
if the time step $k$ decreases by a factor $\ap_k$, then the time spent by both methods increases 
by a factor of $\ap_k$. We remark here that our motivation is to compare the  speed of two methods and thus
their computer programs are run on the same computer.  We understand that the computing time presented
in Tables \ref{T1}--\ref{T2} can be shorten if one uses an advanced computer or does parallel computations,
which however is not our interest here.

\begin{figure}[h!]
\centerline{
\includegraphics[height=6.46cm,width=8.96cm]{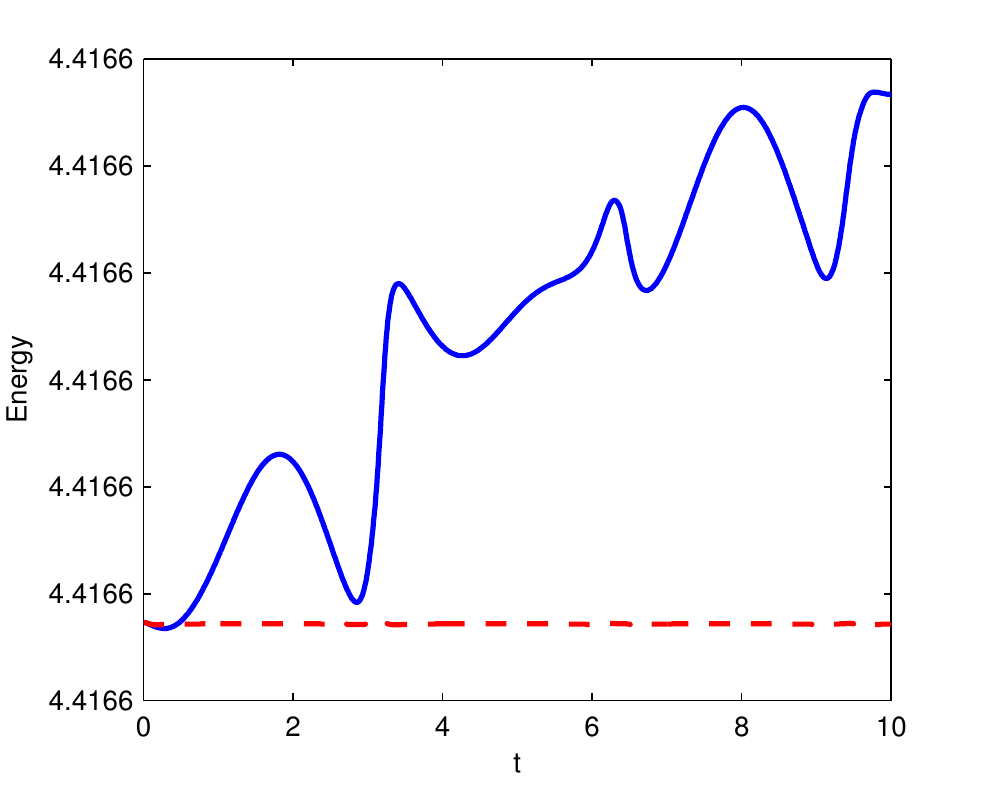}\hspace{-0.2in}
\includegraphics[height=6.46cm,width=8.96cm]{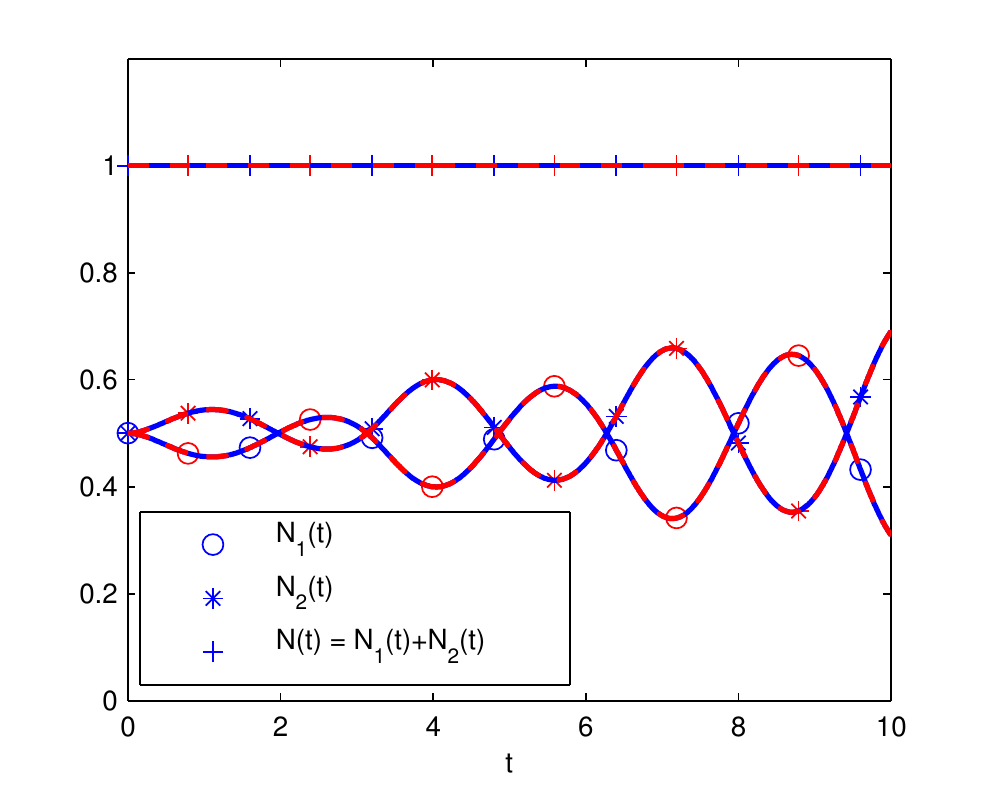}
}
\caption{Dynamics of the energy and the total mass for $t = [0,10]$, where
the mesh size $h = 1/16$ and time step $k = 0.0001$. Solid blue line:  computed from TSADI method;  
dash red line: obtained from our method.}\label{F0}
\end{figure}
In addition,  we study the conservation of  the energy and total mass.  Figure \ref{F0}
shows the time evolution of the energy and total mass for time $t \in[0, 10]$,  where
the mesh size $h =  1/16$ and time step $k = 0.0001$. It shows that both the TSADI and our 
methods conserve the total mass and energy in the discrete level, but our method has a better 
conservation in energy (c.f. Fig. \ref{F0} left).

To further test our method, in Sec. \ref{section5-2}--\ref{section5-3} we will apply it   to study the
dynamical properties of rotating two-component BECs, e.g.,  dynamics of mass, angular momentum
expectation and condensate widths.
Our numerical results will be compared with those reported in \cite{Zhang2007}. In \cite{Zhang2007},
a numerical method was proposed for simulating dynamics of rotating two-component BECs,  in which
the polar  coordinates or cylindrical coordinates were used to resolve the difficulty caused by 
the angular rotational term and a second- or fourth-order finite difference/element discretion is used in the radial direction. Thus, it has low-order accuracy in the radial direction.

\subsection{Dynamics of the mass}
\label{section5-2}

We study dynamics of the mass of each component, i.e., $N_j(t) = \|\psi_j(\cdot, t)\|^2$ for $j = 1, 2$,
and also the total mass $N(t) = N_1(t)+N_2(t)$.  In our simulations,  we solve
two-dimensional CGPEs with the following parameters:
$\lambda = 1$, $\Og = 0.6$  and two-dimensional harmonic potentials are considered with
$\gm_{x, j} = \gm_{y, j} = 1$ ($j = 1, 2$). The initial conditions are chosen as
\beas
\psi_1^0(\bx) = \fl{x+iy}{\sqrt{\pi}} \exp\left(-\fl{x^2+y^2}{2}\right), \qquad
\psi_2^0(\bx) = 0, \qquad \ \bx\in{\mathbb R}^2.
\eeas
That is, initially all atoms are in the first component.
Then we study the dynamics with respect to the following two sets of paramters:
\bea\mbox{(i) \ }
\left(\begin{array}{cc}
 \bt_{11} & \bt_{12} \\
 \bt_{21} & \bt_{22}\end{array}\right) = 500 \left(\begin{array}{cc}
 1.0 & 1.0 \\
 1.0 & 1.0 \end{array}\right); \qquad\mbox{(ii) \ }
 \left(\begin{array}{cc}
 \bt_{11} & \bt_{12} \\
 \bt_{21} & \bt_{22}\end{array}\right) = 500 \left(\begin{array}{cc}
 1.0 & 0.6 \\
 0.6 & 0.8 \end{array}\right).
 \eea
This is one example studied in \cite{Zhang2007}. We remark here that our goal is to test the performance
of our method by comparing it with the available method, and thus we use the same example as that in
\cite{Zhang2007} for the purpose of easy comparison.

\begin{figure}[h!]
\centerline{
(a)\includegraphics[height=5.896cm,width=7.86cm]{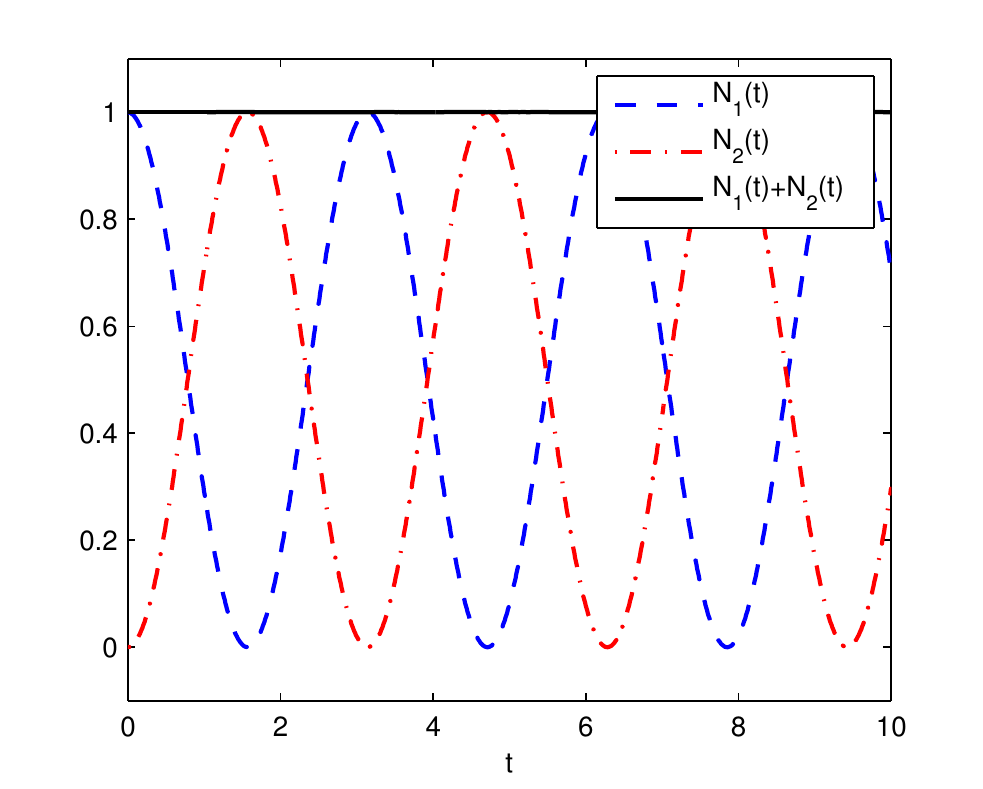}
(b)\includegraphics[height=5.896cm,width=7.86cm]{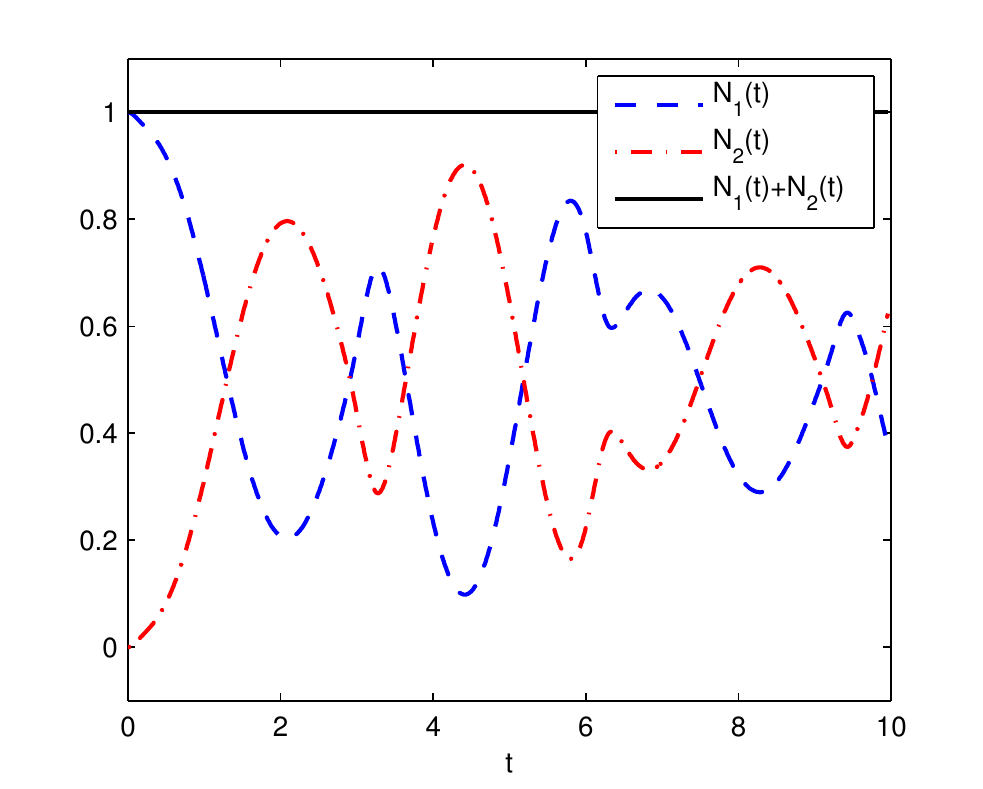}
}
\caption{Time evolution of the mass  $N_j(t) = \|\psi_j(\cdot, t)\|^2$ ($j = 1, 2$) and $N(t) = N_1(t)
+ N_2(t)$ for two sets of interactions parameters: (a) $\bt_{11} = \bt_{12} = \bt_{22} = 500$; (b)
$\bt_{11}$ = 500, $\bt_{22} = 400$, $\bt_{12} = \bt_{21} = 300$.  }\label{Fdensity}
\end{figure}

In our simulations, the computational domain is chosen as ${\mathcal D} = [-8, 8]^2$.  We use mesh size
$h_{\tx} = h_{\ty} = \fl{1}{32}$ and time step $\Dt t = 0.0001$.
Figure \ref{Fdensity} shows the time evolution of  $N_1(t)$, $N_2(t)$ and $N(t)$ for time $t\in[0, 10]$.
From it, we see that when $\bt_{11} = \bt_{12} = \bt_{22}$ (c.f. Fig. \ref{Fdensity}(a)),  the two components
exchange their mass periodically with period $T = \pi/\lambda = \pi$. While when  $\bt_{11} \neq\bt_{12}
\neq \bt_{22}$,  $N_1(t)$ and $N_2(t)$ are not periodical functions (c.f. Fig. \ref{Fdensity}(b)).
In both cases,  the total mass $N(t) = N_1(t) +N_2(t)$ is always conserved. The above
observations are consistent with  the analytical results reported in \cite{Zhang2007, Bao2008}.  Moreover, our
numerical results in Fig. \ref{Fdensity} are the same as those obtained in
\cite{Zhang2007}\footnote{See Figure 3 in \cite{Zhang2007}.} where a numerical method based on
Eulerian coordinates was used. However, our extensive simulations show that
the computing time taken by our method is much shorter than that by the method in \cite{Zhang2007}, 
if the same accuracy is required.

\subsection{Dynamics of angular momentum expectation and condensate widths}
\label{section5-3}

There are  two important quantities in describing the dynamics of rotating two-component BECs:
angular momentum expectation and condensate widths.
In the following, we numerically study their dynamics by applying our method in Sec. \ref{section3}. For
convenience of the readers, we first review the definition of these two quantities; see more information
in \cite{Bao2006, Zhang2007, Wang2007, Bao2008, Bao2013}.

The total angular momentum expectation of two-component BECs is defined  as
\bea
\langle L_z \rangle(t) = \sum_{j = 1}^2 \int_{\mathbb R^d}\psi_j^*(\bx, t) L_z\psi_j(\bx, t)\,
d\bx, 
\qquad t\geq0,
\eea
and the angular momentum expectation of the $j$th component is
\bea
\langle L_z\rangle_j(t) = \fl{1}{N_j(t)} \int_{\mathbb R^d}\psi_j^*(\bx, t) L_z\psi_j(\bx, t)\,d\bx, 
\qquad
t\geq 0, \qquad
\eea
for $j  = 1, 2$.  Usually, the angular
momentum expectation can be used to measure the vortex flux.  The condensate width of
two-component BECs in $\ap$-direction ($\ap = x, y$ or $z$) is defined as
\bea
\sigma_\ap = \sqrt{\dt_{\ap}(t)} = \sqrt{\dt_{\ap, 1}(t) + \dt_{\ap, 2}(t)}, \quad \  t\geq 0, \qquad
\ap = x, \  y\  \mbox{or}\  z,
\eea
where
\bea
\dt_{\ap, j}(t) = \langle \ap^2\rangle_j(t) = \int_{{\mathbb R}^d} \ap^2|\psi_j(\bx, t)|^2 d\bx, \quad \
t\geq 0, \qquad j = 1, 2.
\eea

To study the dynamics of angular momentum expectation and condensate widths, we choose the
following parameters in the CGPEs (\ref{DCGPEs1})--(\ref{DCGPEs2}):
 $d = 2$,   $\Omega = 0.6$, $\lambda = 1$ and
\bea
 \left(\begin{array}{cc}
 \bt_{11} & \bt_{12} \\
 \bt_{21} & \bt_{22}\end{array}\right) = 400 \left(\begin{array}{cc}
 1.0 & 0.97 \\
 0.97 & 0.94 \end{array}\right).
 \eea
The initial conditions are taken as
\beas
\psi_j^0(\bx) = \fl{x+iy}{\sqrt{2\pi}}\exp\left(-\fl{x^2+y^2}{2}\right), \qquad \bx\in{\mathbb R}^2,
\quad j = 1, 2.
\eeas
The computational domain is chosen as ${\mathcal D} = [-24,24]^2$ with the mesh
size $h_{\tx} = h_{\ty} = \fl{3}{64}$ and time step is $k = 0.0001$.

\begin{figure}[h!]
\centerline{
(a)\includegraphics[height=5.896cm,width=7.86cm]{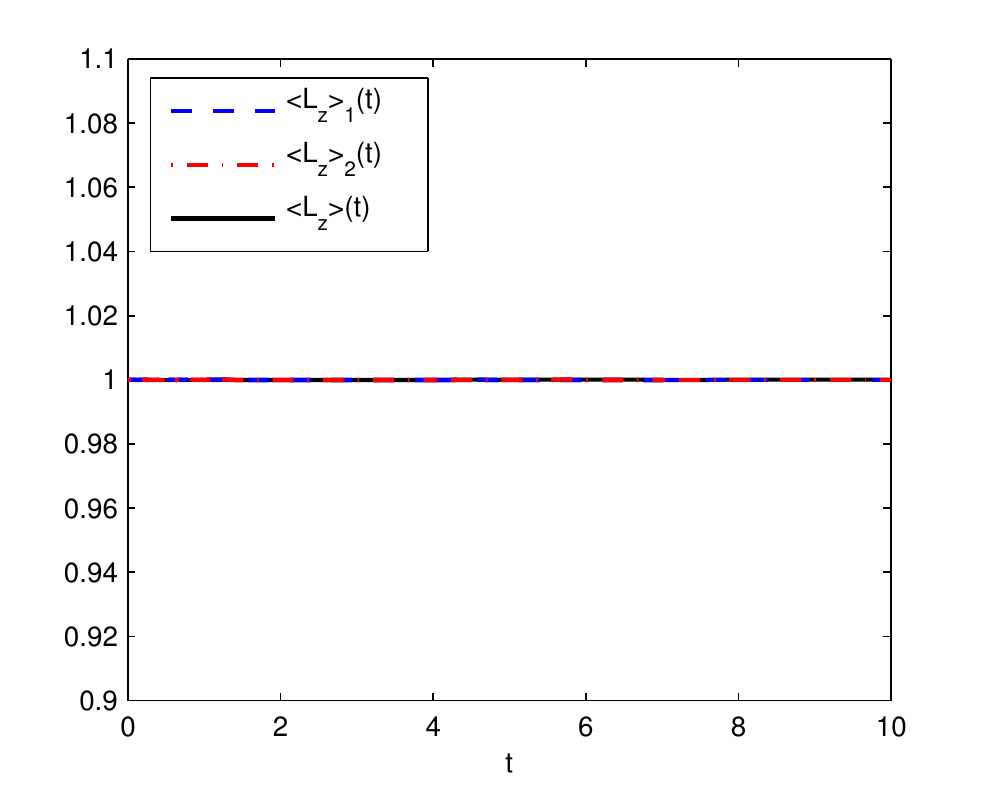}
(b)\includegraphics[height=5.896cm,width=7.86cm]{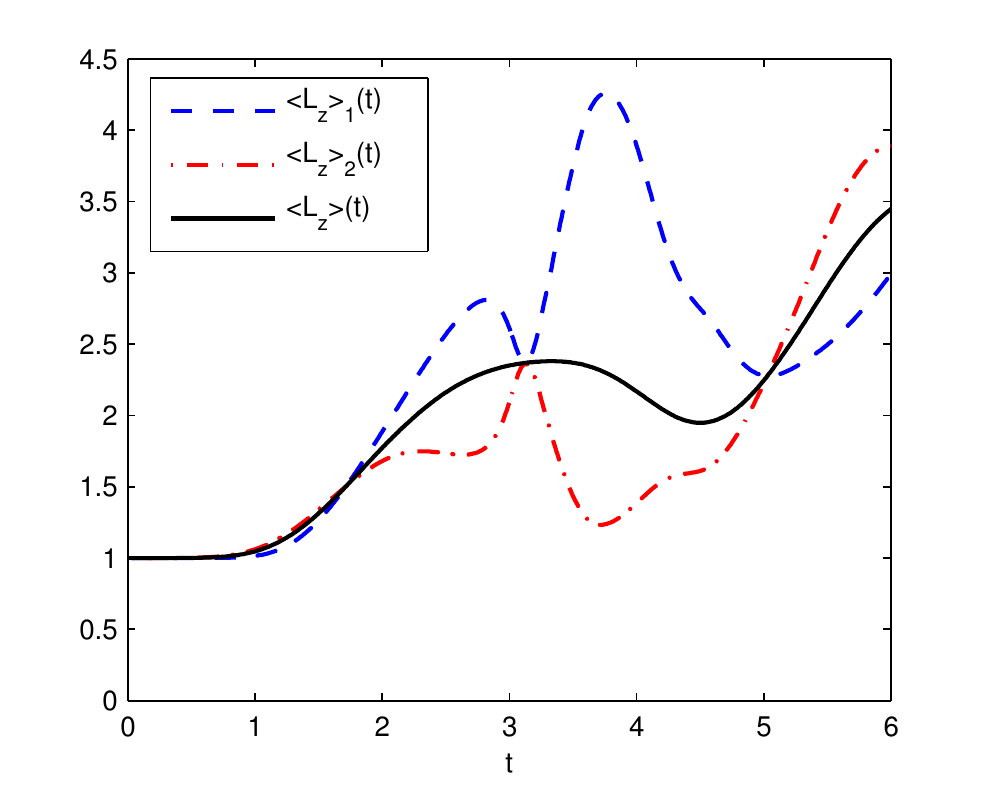}
}
\caption{Time evolution of angular momentum expectation  $\langle L_z\rangle_j(t)$ ($j = 1, 2$) and
$\langle L_z\rangle$(t) for two sets of trapping frequencies:  (a) $\gm_{x, j} = \gm_{y,j} = 1$ ($j = 1, 2$); 
(b) $\gm_{x, 1} = \gm_{y,1} = 1$, $\gm_{x, 2} = 1.05, \gm_{y,2} = 0.9$.}\label{Flz}
\end{figure}

Figure \ref{Flz} presents the dynamics of angular momentum expectations for two sets of
trapping frequencies:  (i) $\gm_{x, 1} = \gm_{y,1} = \gm_{x, 2} = \gm_{y,2} = 1$; (ii)
$\gm_{x, 1} = \gm_{y,1} = 1$, $\gm_{x, 2} = 1.05, \gm_{y,2} = 0.9$.  It shows that the total angular
momentum expectation $\langle L_z\rangle (t)$ is conserved as long as both external trapping
potentials $V_1(\bx)$ and $V_2(\bx)$ in (\ref{potential}) are symmetric (c.f. Fig. \ref{Flz}(a)).  While
when $\lambda \neq 0$,  if at least one 
of the external potentials is asymmetric,  none of $\langle L_z\rangle_1(t)$,
$\langle L_z\rangle_2(t)$ and $\langle L_z\rangle(t)$ is conserved  (c.f.
Fig. \ref{Flz}(b)).   In addition,  the results in Fig. \ref{Flz} are the same  as those reported in
\cite{Zhang2007}\footnote{See Figure 4  in \cite{Zhang2007}.} but the computing time used 
by our method is much less.

\begin{figure}[h!]
\centerline{
(a)\includegraphics[height=5.996cm,width=7.96cm]{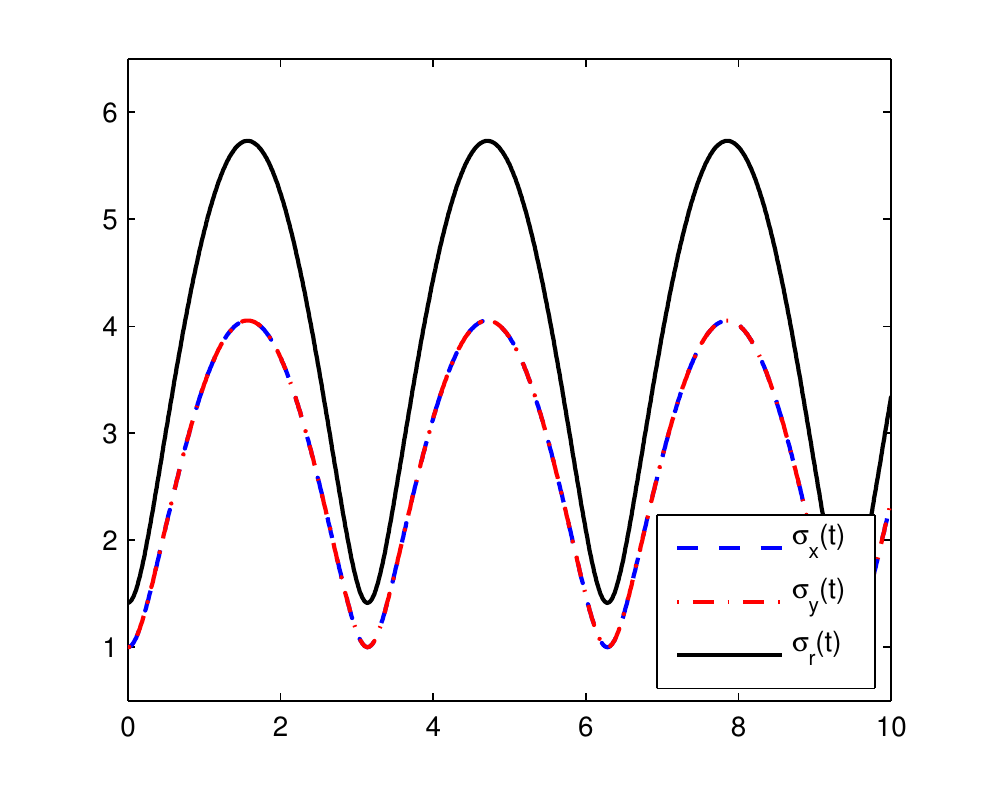}
(b)\includegraphics[height=5.996cm,width=7.96cm]{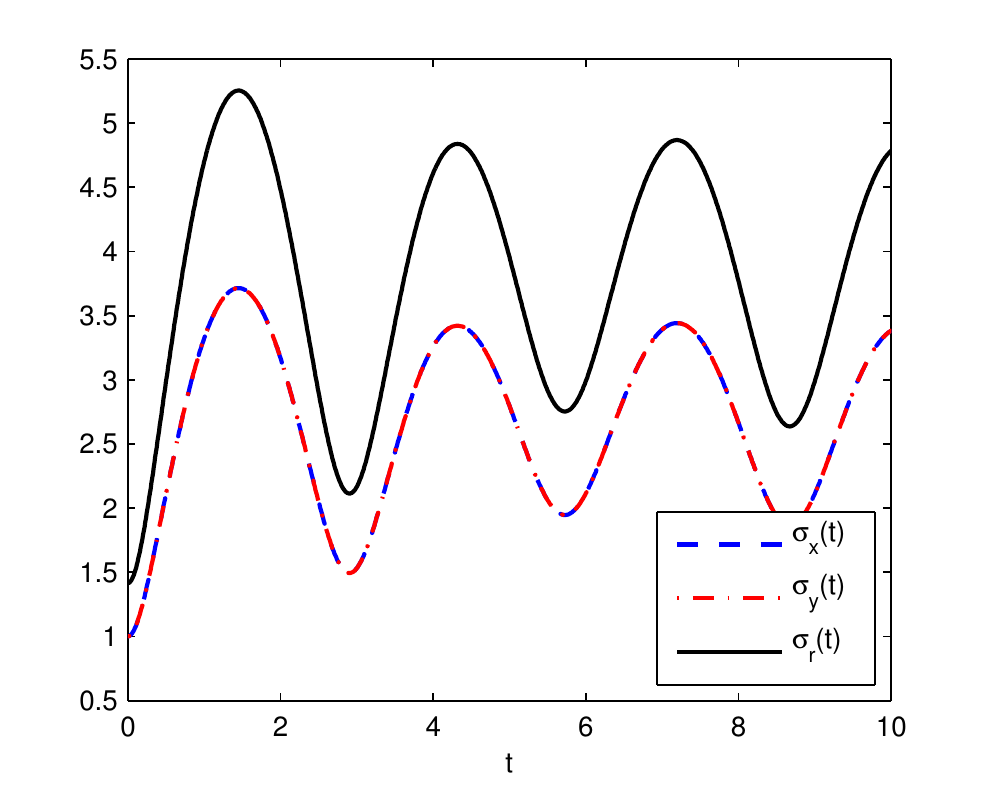}
}
\caption{Time evolution of condensate widths $\sigma_x(t)$, $\sigma_y(t)$ and $\sigma_r(t)$ for two
sets of trapping frequencies:  (a) $\gm_{x, j} = \gm_{y,j} = 1$ ($j = 1, 2$); (b) $\gm_{x, 1} = \gm_{y,1} = 1$, $\gm_{x, 2} = \gm_{y,2} = 1.2$.}\label{Fwidth}
\end{figure}

Figure \ref{Fwidth} shows the dynamics of condensate widths for $\sigma_x(t)$, $\sigma_y(t)$ and
$\sigma_r(t) := \sqrt{\sigma_x^2(t) + \sigma^2_y(t)}$ for two sets of
trapping frequencies:  (i) $\gm_{x, 1} = \gm_{y,1} = \gm_{x, 2} = \gm_{y,2} :=\gamma = 1$; (ii)
$\gm_{x, 1} = \gm_{y,1} = 1$, $\gm_{x, 2} = \gm_{y,2} = 1.2$.  From it, we see that when
the two components have the same external trapping potentials, the condensate widths
$\sigma_x(t)$, $\sigma_y(t)$ and $\sigma_r(t)$ are periodic functions with period
$T = \pi/\gamma = \pi$ (c.f. Fig. \ref{Fwidth}(a)). If the potential $V_1(\bx)\neq V_2(\bx)$ in (\ref{potential}),
the condensate widths are not periodic functions. The above results are consistent with those in
\cite{Zhang2007}\footnote{See Figure 5 in \cite{Zhang2007}.}.

\subsection{Dynamics of vortex lattices}

In this section, we apply our method to study the dynamics of  vortex lattices in rotating
two-component BECs. The initial data are taken as the stationary vortex lattices, which are
computed by choosing $\Og = 0.9$,  $\lambda = 0$,  $\gm_{x, j} = \gm_{y, j} = 1$ ($j = 1, 2$) and
\beas
\left(\begin{array}{c c}
 \bt_{11} & \bt_{12} \\
 \bt_{21} & \bt_{22}\end{array}\right) = 500\left(\begin{array}{cc}
 1.0 & -0.25 \\
-0.25 & 1.0 \end{array}\right).
\eeas
Due to the attractive interaction between two components, i.e., $\bt_{12} = \bt_{21} < 0$,  
initially  the stationary vortex lattices are exactly the same.  Then at time $t = 0$, 
\begin{enumerate}
\item[] Case (i).  Change the symmetric external potentials to asymmetric by setting
$\gamma_{x,1} = \gamma_{y,2} = 1.05$ and $\gamma_{x,2} = \gamma_{y,1} = 0.95$;
\item[] Case (ii).   Turn on the external driving field by setting $\lambda = 1$. 
\end{enumerate}
Then we study the dynamics of vortex lattices. The computational domain is chosen as ${\mathcal D}
= [-24, 24]^2$ with $h_{\tx} = h_{\ty} = \fl{3}{32}$ and the time step $k = 0.0001$.

\begin{figure}[h!]
\centerline{
\includegraphics[height=3.86cm,width=4cm]{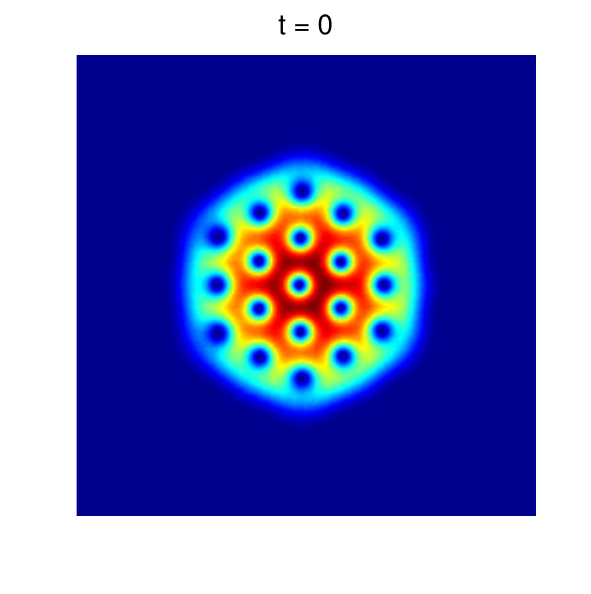}\hspace{-0.4cm}
\includegraphics[height=3.86cm,width=4cm]{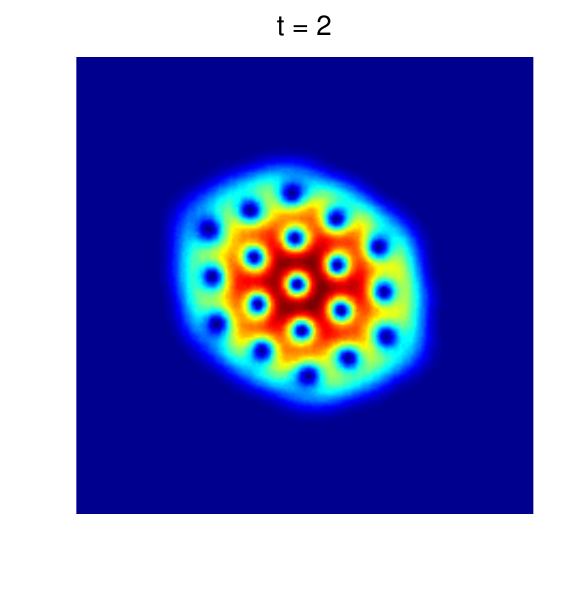}\hspace{-0.4cm}
\includegraphics[height=3.86cm,width=4cm]{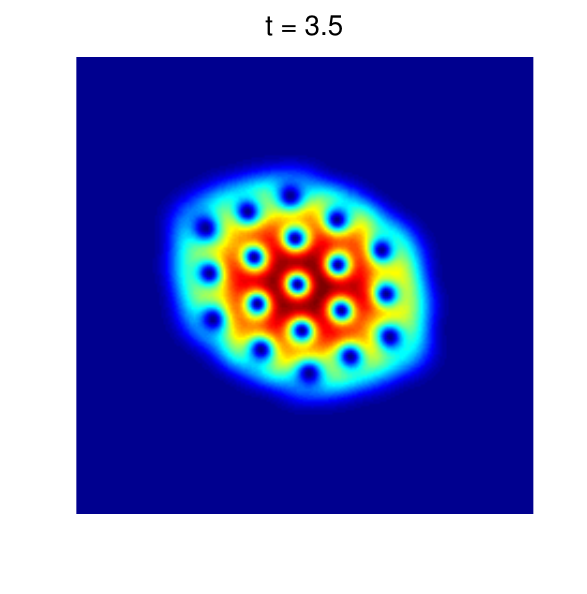}\hspace{-0.4cm}
\includegraphics[height=3.86cm,width=4cm]{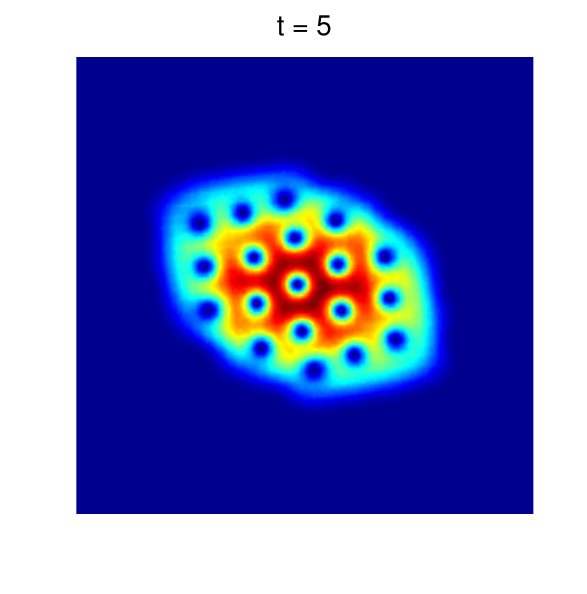}
}
\centerline{
\includegraphics[height=3.86cm,width=4cm]{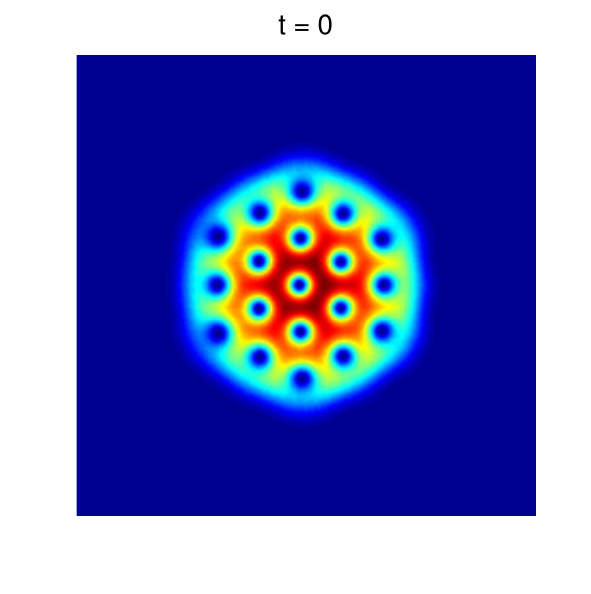}\hspace{-0.4cm}
\includegraphics[height=3.86cm,width=4cm]{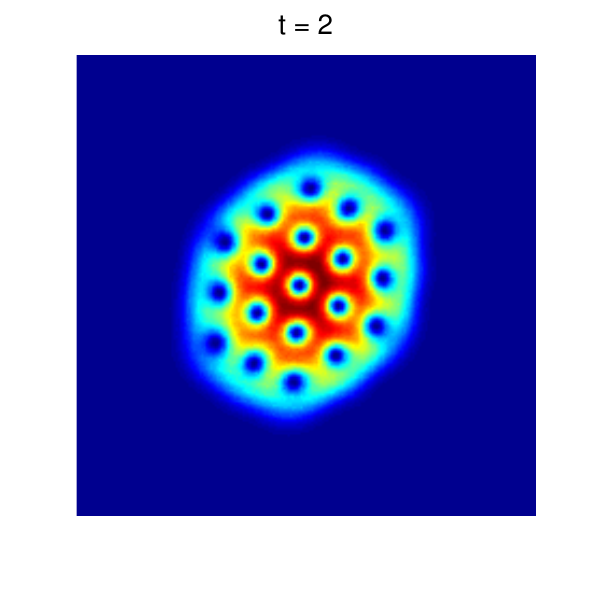}\hspace{-0.4cm}
\includegraphics[height=3.86cm,width=4cm]{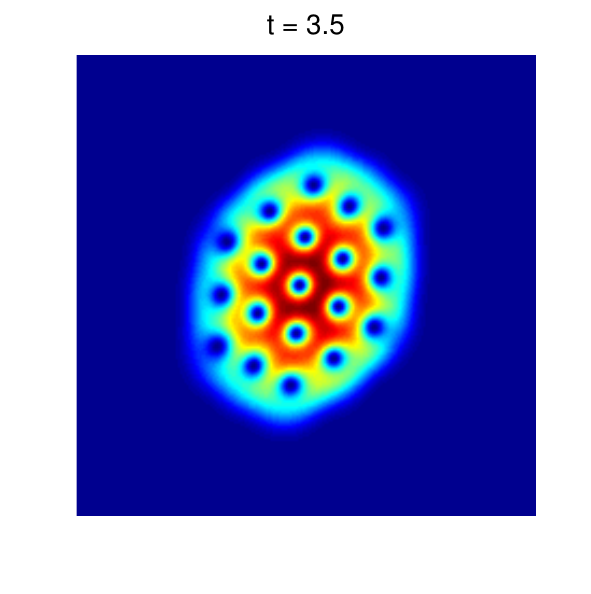}\hspace{-0.4cm}
\includegraphics[height=3.86cm,width=4cm]{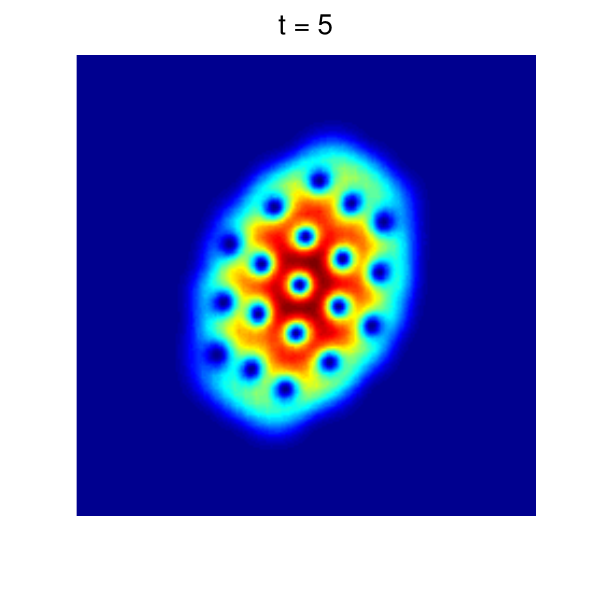}
}
\caption{Contour plots of the density $|\psi_1|^2$ (top row) and $|\psi_2|^2$ (bottom row)
at different time for  Case (i). Displayed domain $(x,y) \in [-10, 10]^2$.}\label{Flattice1}
\end{figure}
\begin{figure}[h!]
\centerline{
\includegraphics[height=3.86cm,width=4cm]{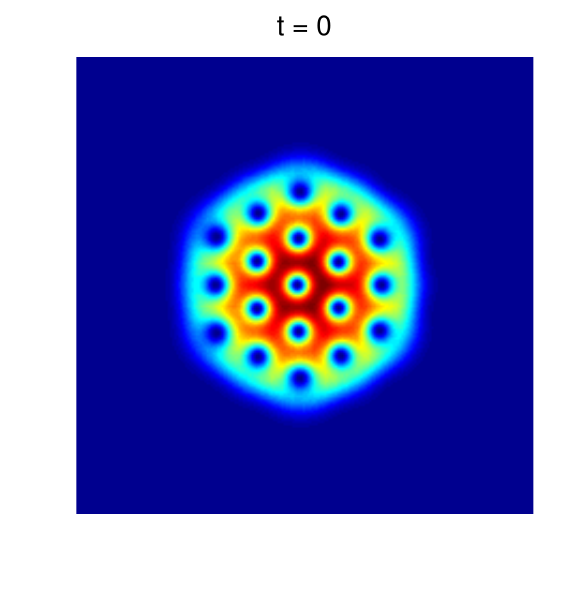}\hspace{-0.4cm}
\includegraphics[height=3.86cm,width=4cm]{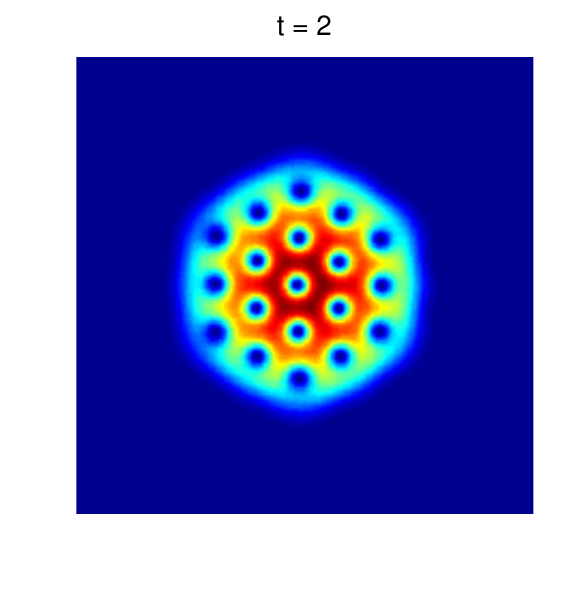}\hspace{-0.4cm}
\includegraphics[height=3.86cm,width=4cm]{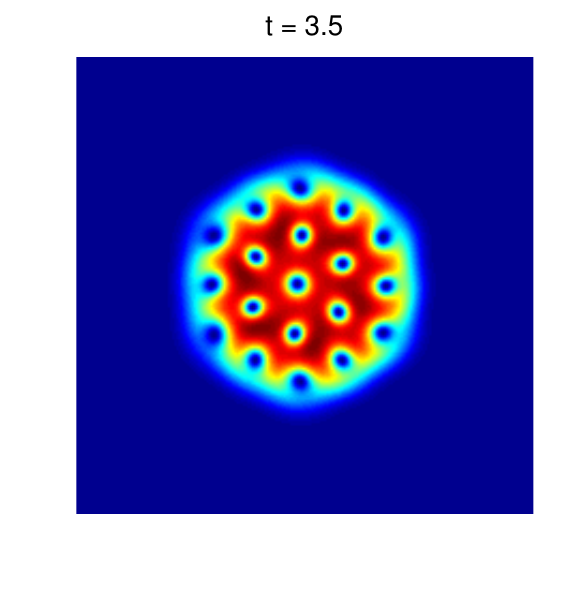}\hspace{-0.4cm}
\includegraphics[height=3.86cm,width=4cm]{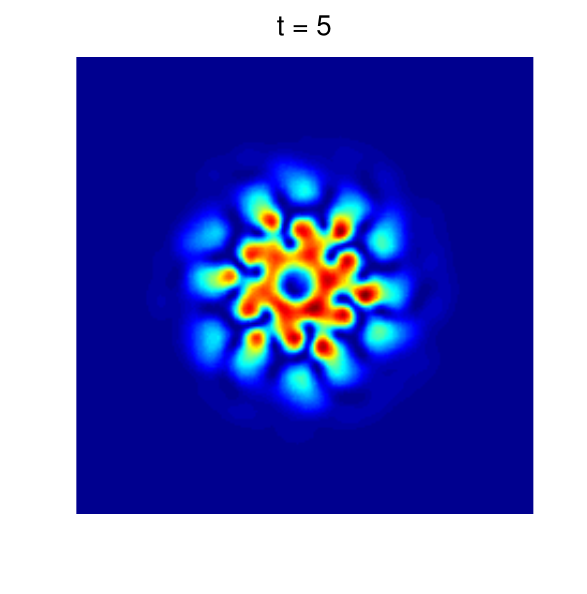}
}
\centerline{
\includegraphics[height=3.86cm,width=4cm]{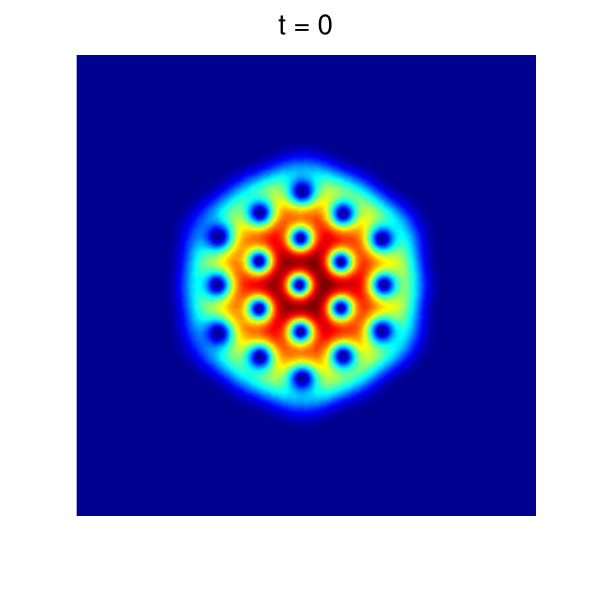}\hspace{-0.4cm}
\includegraphics[height=3.86cm,width=4cm]{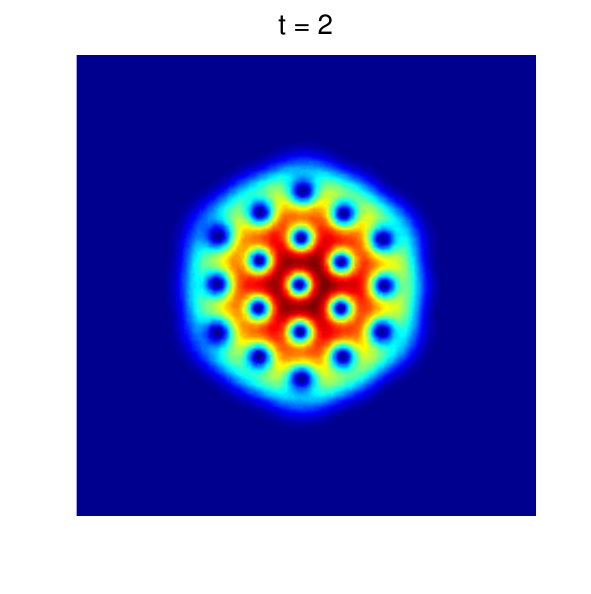}\hspace{-0.4cm}
\includegraphics[height=3.86cm,width=4cm]{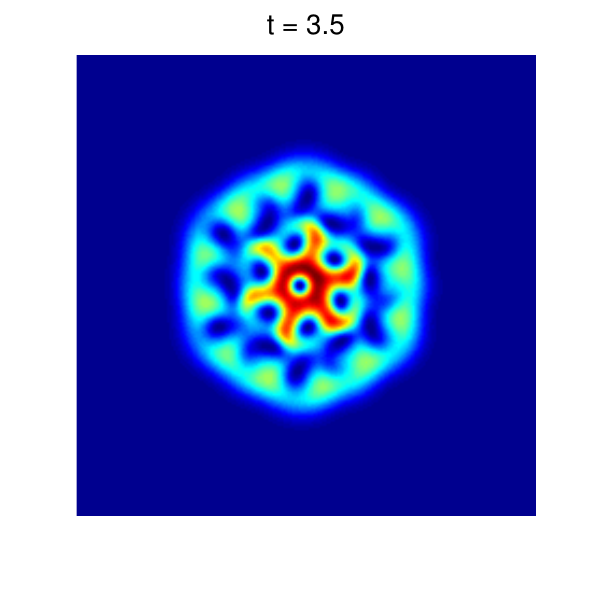}\hspace{-0.4cm}
\includegraphics[height=3.86cm,width=4cm]{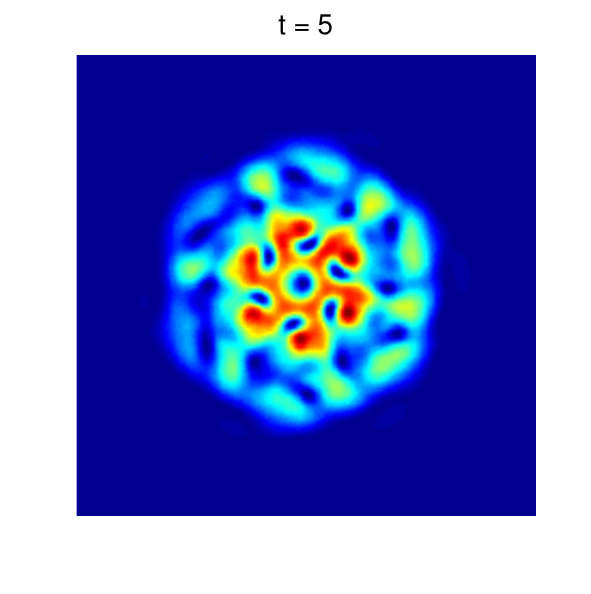}
}
\caption{Contour plots of the density $|\psi_1|^2$ (top row) and $|\psi_2|^2$ (bottom row)
at different time for  Case (ii). Displayed domain $(x,y) \in [-10, 10]^2$.}\label{Flattice2}
\end{figure}

Figures \ref{Flattice1}--\ref{Flattice2} show the contour plots of the density $|\psi_1|^2$ and
$|\psi_2|^2$ at different time $t$ in Case (i) and (ii), respectively, where the displayed domain is
$(x, y) \in[-10, 10]^2$. At $t = 0$, the vortex lattices are identical and there are
19 vortices in each lattice.  In  Case (i),  the lattices rotate periodically due to the anisotropy of
the external potentials and the number of vortices is conserved during the dynamics (c.f.  Fig.
\ref{Flattice1}). While in Fig. \ref{Flattice2}, we see that the external driving field eventually
destroys the pattern of vortex lattices.

\section{Summary}
\label{section6}

We proposed an  efficient numerical method to solve the coupled Gross--Pitaevskii equations
(CGPEs) with both angular momentum rotation term and external driving field term, which well
describes the dynamics of rotating two-component Bose--Einstein condensations (BECs) with an 
internal  Josephson junction.  We introduced a rotating Lagrangian coordinate transformation and then 
eliminated the angular momentum  rotation term in the CGPEs. Under the new coordinates, we 
proposed  a time-splitting sine pseudospectral method to simulate the dynamics of rotating 
two-component BECs. To efficiently treat the external driving field term, we applied  
a linear transformation so that it does not cause any  extra computational complexity. Compared to the
methods  in literature has higher order spatial accuracy but  requires less memory cost and 
computational cost.  It can be easily implemented in practice. We then numerically examined
the conservation of the angular momentum expectation  and studied the dynamics of condensate
widths and center of mass for different angular velocities. In addition, the dynamics of vortex lattice
in rotating two-component BECs were investigated. Numerical studies showed that our method is 
very effective.

\vskip 20pt
\noindent{\large\bf Acknowledgements \ }  \ Q. Tang and Y. Zhang would like to express their sincere thanks to
Prof. Weizhu Bao for the fruitful discussion on the project. This work was partially supported by
the Singapore A*STAR SERC Grant No. 1224504056 (Q. Tang) and
by the Simons Foundation Award No. 210138 (Y. Zhang).

\bigskip
\bigskip

\noindent{\large{\bf References}}

\end{document}